\documentclass[11pt]{amsart}
\usepackage[margin=2.6cm]{geometry}
\usepackage{amsmath,amssymb,amsthm,mathtools}
\usepackage{enumerate}
\usepackage{xcolor}
\usepackage{array}
\usepackage{pgfplots}
\pgfplotsset{compat=1.16}
\usetikzlibrary{patterns}
\definecolor{lava}{rgb}{0.81, 0.06, 0.13}
\usepackage[colorlinks=true,linkcolor=lava,citecolor=blue,urlcolor=blue!60!black]{hyperref}

\numberwithin{equation}{section}

\newtheorem{theorem}{Theorem}[section]
\newtheorem{proposition}{Proposition}[section]
\newtheorem{lemma}{Lemma}[section]
\newtheorem{corollary}{Corollary}[section]
\newtheorem{conjecture}{Conjecture}[section]

\theoremstyle{definition}
\newtheorem{definition}{Definition}[section]

\theoremstyle{remark}
\newtheorem{remark}{Remark}[section]

\newcommand{\R}{\mathbb{R}}

\newcommand{\Ha}{H^1(\R^N)}
\newcommand{\Lt}{L^2(\R^N)}
\newcommand{\e}{\mathrm{e}}
\newcommand{\dd}{\,\mathrm{d}}

\newcommand{\Fou}{\mathcal{F}}
\newcommand{\cS}{\mathcal{S}}
\newcommand{\Lsub}{\underline{a}}
\newcommand{\da}{\mathfrak{a}}
\newcommand{\Mo}{\mathcal{M}}
\newcommand{\Do}{\mathcal{D}}

\begin{document}

\title[Scattering thresholds for NLS with decaying damping]{Scattering thresholds for the NLS equation with polynomially decaying damping}

\author{Makram Hamouda}
\address{(M. Hamouda) Department of Basic Sciences, Deanship of Preparatory Year and Supporting Studies, Imam Abdulrahman Bin Faisal University, P.~O.~Box 1982, Dammam, Saudi Arabia.}
\email{mmhamouda@iau.edu.sa}
\urladdr{https://orcid.org/0000-0003-1598-3962}
\author{Mohamed Majdoub}
\address{(M. Majdoub) Department of Mathematics, College of Science, Imam Abdulrahman Bin Faisal University, P.~O.~Box 1982, Dammam, Saudi Arabia.}
\address{Basic and Applied Scientific Research Center, Imam Abdulrahman Bin Faisal University, P.~O.~Box 1982, 31441, Dammam, Saudi Arabia.}
\email{mmajdoub@iau.edu.sa}
\email{med.majdoub@gmail.com}
\email{mohamed.majdoub@fst.rnu.tn}
\urladdr{https://orcid.org/0000-0001-6038-1069}

\begin{abstract}
We study the focusing nonlinear Schr\"odinger equation with a linear damping of strength $a(1+t)^{-\theta}$, where $a>0$ and $\theta\ge0$, in the whole energy-subcritical range. For $\theta>0$, the time-averaged damping vanishes, so that the scattering theory developed in our previous work does not apply. We show that the long-time dynamics is governed by the time integrability of an effective coupling, which decays faster than any power when $\theta<1$, polynomially when $\theta=1$, and not at all when $\theta>1$. In the mass-subcritical case, all solutions are global. They scatter when $\theta<1$ and, for $\theta=1$, above an explicit exponent depending on $a$. On the other hand, no non-trivial solution is asymptotically free below a second explicit exponent when $\theta=1$, and below the short-range exponent $1+2/N$ when $\theta>1$. In the mass-critical case, solutions are global below the ground state. Moreover, for $\theta\le1$, every solution which survives until its decaying mass becomes subcritical is global and scatters. In the intercritical case, we prove small data scattering and, for $\theta\le1$, a dichotomy between scattering and grow-up of the gradient. Finally, in the scale-invariant case $\theta=1$, we analyse the range between the two thresholds. We reduce the equation exactly to a nonlinear Schr\"odinger equation with a decaying repulsive harmonic potential, and we construct approximate solitary waves with an integrable error. This leads to a conjecture, supported by numerical simulations, on the sharpness of the scattering threshold.
\end{abstract}

\maketitle

\section{Introduction and main results}\label{S1}

This paper is devoted to the Cauchy problem for the focusing nonlinear Schr\"odinger equation with a linear damping whose strength decays polynomially in time,
\begin{equation}\label{Main-New}
\begin{cases}
i\partial_t u+\Delta u+ i\,\da(t)\,u=-|u|^{p-1}u, & (t,x)\in(0,\infty)\times\R^N,\\[2pt]
u(0,x)=u_0(x), & x\in\R^N,
\end{cases}
\qquad \da(t)=\dfrac{a}{(1+t)^{\theta}},
\end{equation}
where $N\ge1$, $a>0$, $\theta\ge0$ and $1<p<1+\frac{4}{(N-2)_+}$. Here and below, $x_+=\max(x,0)$ for $x\in\R$, with the convention $0^{-1}=\infty$. Throughout the paper, the damping \emph{function} is denoted by the Fraktur letter $\da$, while $a$ always denotes the positive \emph{constant} that measures its strength in \eqref{Main-New}. The term $i\da(t)u$ models a loss mechanism whose intensity varies with time. Equations of this type arise in nonlinear optics, where the time-dependent damping corresponds to an absorption coefficient. They also arise in the theory of Bose--Einstein condensates, where $\da(t)$ is a loss coefficient describing the escape of atoms from the condensate; see \cite{Fibich2001,FengZhaoSun2014,PerezGarcia1995,BamriTayachi2023}, \cite{KivsharMalomed1989} and the references therein.

The model \eqref{Main-New} was singled out in \cite[Remark 2.6 and (2.10)]{HM-MJM2025}, which deals with the general equation
\begin{equation}\label{general}
i\partial_t u+\Delta u+ i\,\da(t)\,u=-|u|^{p-1}u,\qquad \da\in C([0,\infty);[0,\infty)).
\end{equation}
We set
\begin{equation}\label{A}
A(t)=\int_0^t \da(s)\dd s,\qquad \Lsub=\inf_{t>0}\frac{A(t)}{t}.
\end{equation}
All the scattering results of \cite{HM-MJM2025} were obtained under the assumption $\Lsub>0$, which forces the damping to be, on average, at least as strong as a positive constant. For the damping in \eqref{Main-New}, a direct integration gives
\begin{equation}\label{A-explicit}
A(t)=\begin{cases}
\dfrac{a}{1-\theta}\big((1+t)^{1-\theta}-1\big), & \theta\neq1,\\[8pt]
a\log(1+t), & \theta=1.
\end{cases}
\end{equation}
We infer that $\Lsub=a$ when $\theta=0$, while $\Lsub=0$ for every $\theta>0$. The case $\theta>0$ is therefore not covered by \cite{HM-MJM2025}, and its study was announced there as a natural continuation. The purpose of the present work is to give a fairly complete answer. We establish global existence, scattering and non-scattering results for \eqref{Main-New} in the whole energy-subcritical range, and we analyse in detail the only range that our methods leave open.

\subsection{Background}
When $\da\equiv0$, \eqref{general} is the classical focusing NLS equation. Its local and global theory in the energy space, as well as its scattering theory, are well understood; see the monographs \cite{Cazenave2003,Tao2006,LinaresPonce2015,BahouriCheminDanchin2011} and the papers \cite{GinibreVelo1979,GinibreVelo1985,Kato1987,Tsutsumi1985,CazenaveWeissler1990,NakanishiOzawa2002}. Three exponents organise this theory. For $1<p<1+4/N$ (the mass-subcritical case), all $H^1$ solutions are global and bounded in $H^1$. Non-trivial solutions are never asymptotically free when $p\le 1+2/N$. This was proved by Strauss \cite{Strauss1974} for $1<p\le1+2/N$ when $N\ge2$ and for $1<p\le2$ when $N=1$, and by Barab \cite{Barab1984} in the remaining one-dimensional range $2<p\le3$; see also \cite{Glassey1973} for the analogous result for the wave equation, and \cite{TsutsumiYajima1984}, \cite[Section 7.5]{Cazenave2003}. At the borderline $p=1+2/N$, small solutions exhibit a logarithmic phase correction, that is, a modified scattering \cite{OzawaLR1991,GinibreOzawa1993,HayashiNaumkin1998}. In the focusing case, the solitary waves $\e^{it}Q(x)$ show that scattering fails for large data in the whole mass-subcritical range. These solitary waves are orbitally stable \cite{CazenaveLions1982,Weinstein1985,Weinstein1986}. In the defocusing mass-subcritical case, scattering for data in the weighted space $\Sigma$ holds above the Strauss exponent \cite{Strauss1981,Tsutsumi1985,CazenaveWeissler1992}, while scattering for general $H^1$ data is a well-known open problem; see \cite{KMMV2017} for conditional results. At the mass-critical exponent $p=1+4/N$, solutions with $\|u_0\|_{L^2}<\|Q\|_{L^2}$, where $Q$ is the ground state, are global \cite{Weinstein1983} and scatter \cite{Dodson2015}. In the intercritical range $1+4/N<p<1+4/(N-2)_+$, small data scatter, and the sharp threshold for scattering is described in terms of the ground state \cite{HolmerRoudenko2008,DuyckaertsHolmerRoudenko2008,FangXieCazenave2011,Guevara2014,KenigMerle2006}; see also \cite{DinhMajdoubSaanouni2023} for an inhomogeneous variant. Finite time blow-up may occur for large data \cite{Glassey1977}.

For a constant damping $\da(t)\equiv a>0$, the equation was first studied by M.~Tsutsumi \cite{TsutsumiM1984,TsutsumiM1990}, who proved that blow-up still occurs when the damping is weak compared with the size of the data. Global existence for a large damping was obtained by Ohta and Todorova \cite{OhtaTodorova2009}; see also \cite{Fibich2001,Darwich2012,Dinh2021}, \cite{Saanouni2020} for refined blow-up and global existence criteria. The asymptotic behaviour for a constant damping was investigated by Inui \cite{Inui2019}, who proved scattering in the energy space, in the natural sense recalled in Definition \ref{Def-Scat} below, in both the focusing and the defocusing regimes; sharp convergence rates were obtained very recently in \cite{TakagiTakizawa2026}. Dissipative nonlinear terms of a different nature were studied by Kita and Shimomura \cite{KitaShimomura2007,KitaShimomura2009}, by Cazenave, Han and Naumkin \cite{CazenaveHanNaumkin2021}, and by Antonelli, Carles and Sparber \cite{AntonelliCarlesSparber2015}.

For genuinely time-dependent dampings, much less is known. The global theory in $H^1$ for a general $\da\in L^1_{\rm loc}$ was investigated by Feng, Zhao and Sun \cite{FengZhaoSun2014}. Darwich \cite{Darwich2018} obtained global existence for positive time-dependent dissipation coefficients. Cavalcanti, Corr\^ea, Domingos Cavalcanti and Tebou \cite{CCDT2017} studied well-posedness and energy decay for the defocusing problem with a damping coefficient which may vanish at infinity. Time-dependent dampings were also considered for Choquard and fourth-order Schr\"odinger equations \cite{FengSaanouni2021,Saanouni2021}; see also \cite{Chergui2021} for the Choquard equation with a constant damping. In \cite{HM-CPAA2023}, we proved that a sufficiently strong time-dependent damping yields global solutions, while a weak one does not prevent blow-up. There, the divergence $\int_0^\infty \da(s)\dd s=\infty$ appears as a natural condition. In \cite{HM-MJM2025}, we developed an energy scattering theory under the condition $\Lsub>0$, in the spirit of \cite{Inui2019}. The influence of an inverse-square potential and of an inhomogeneous nonlinearity was analysed in \cite{HMS-ZAMP2025}. Closely related problems for nonlinear Schr\"odinger equations with exponential nonlinearities, and in the energy-critical regime, were considered in \cite{DinhKeraaniMajdoub2020,CIMM2009}. 

A different and very elegant point of view was adopted by Bamri and Tayachi \cite{BamriTayachi2023} (see also \cite{Bamri2025}). They considered dampings with the asymptotic behaviour $t\,\da(t)\to\gamma\ge0$ and oscillating data. The starting point of \cite{BamriTayachi2023} is the transformation $v=\e^{A(t)}u$, which converts \eqref{general} into the nonlinear Schr\"odinger equation
\begin{equation}\label{main-bis}
i\partial_t v+\Delta v+h(t)|v|^{p-1}v=0,\qquad h(t)=\e^{-(p-1)A(t)},
\end{equation}
with a time-dependent coupling in front of the nonlinearity. Equations of the form \eqref{main-bis} with a general $h$ were previously considered by Fanelli \cite{Fanelli2009}, by Gonda, Machihara and Ozawa \cite{GondaMachiharaOzawa2014}, by Cazenave and Scialom \cite{CazenaveScialom2010} and by Carles \cite{Carles2011}. Working in Strichartz spaces with Lorentz norms in time, Bamri and Tayachi proved global existence and scattering for small data when $N\ge2$ and $h\in L^{1/\delta,\infty}(0,\infty)$ with $0<\delta<1$, in a range of exponents depending on $\delta$ \cite[Theorem 1.2 and Corollaries 1.4, 1.6]{BamriTayachi2023}. Using the pseudo-conformal transformation and the method of Cazenave and Naumkin \cite{CazenaveNaumkin2017}, they also obtained scattering for a special class of non-vanishing oscillating data whenever $p-1>2/(N+\gamma)$ \cite[Theorem 1.9]{BamriTayachi2023}. Finally, they proved that solutions do not scatter when $p-1<2/(N+\gamma)$, for real $\lambda\neq0$, data in the weighted space $\Sigma$ and $p\le2$ \cite[Theorem 1.7]{BamriTayachi2023}, and that this critical value is unchanged under $L^1$-perturbations of the damping. Their normalisation is $i\partial_tu+\Delta u+\frac{i}{2}\da_{\rm BT}(t)u+\lambda|u|^{p-1}u=0$. Consequently, \eqref{Main-New} with $\theta=1$ corresponds to $\da_{\rm BT}(t)=2a/(1+t)$, that is, to $\gamma=2a$ and $\lambda=1$.

It is instructive to compare \eqref{Main-New} with the damped wave equation $\partial_t^2w-\Delta w+\mu(1+t)^{-\theta}\partial_tw=f(w)$. Wirth's classification \cite{Wirth2006,Wirth2007} distinguishes three regimes. In the \emph{effective} range $\theta<1$, the damped wave behaves like a heat equation. In the \emph{non-effective} range $\theta>1$, it behaves like a free wave. In the \emph{scale-invariant} case $\theta=1$, the size of $\mu$ decides between the two behaviours. We refer to \cite{Nishihara2011,LinNishiharaZhai2012,DAbbicco2015} for the corresponding critical exponents, and to \cite{HamoudaHamza2021a,HamoudaHamza2021b} for blow-up results in the scale-invariant case with combined nonlinearities. Our results show that a similar trichotomy governs the Schr\"odinger problem \eqref{Main-New}. However, the mechanism is of a completely different nature, because the Schr\"odinger flow is conservative. The damping does not change the dispersion. Instead, it reduces the effective strength of the nonlinearity through the factor $h(t)$ in \eqref{main-bis}.

\subsection{Notation, setting and basic definitions}\label{S1-setting}

We write $\|\cdot\|_{L^r}$ for the norm of $L^r(\R^N)$, $m=\|u_0\|_{L^2}$, and $(f,g)=\int_{\R^N}f\bar g\dd x$. For $1\le r\le\infty$, $r'$ denotes the conjugate exponent of $r$. The Fourier transform is $\Fou f(\xi)=\hat f(\xi)=(2\pi)^{-N/2}\int\e^{-ix\cdot\xi}f(x)\dd x$. We set $2^*=\frac{2N}{N-2}$ if $N\ge3$, and $2^*=\infty$ if $N=1,2$. We use the norms $\|f\|_{W^{1,r}}=\|f\|_{L^r}+\|\nabla f\|_{L^r}$ and $\|f\|_{H^1}=\|f\|_{W^{1,2}}$, so that the Strichartz constants below apply without additional factors. We denote by $B(x_0,\rho)$ the open ball of $\R^N$ with centre $x_0$ and radius $\rho$, and by $\Sigma=\{f\in\Ha:\ |x|f\in\Lt\}$ the usual weighted space. The constant $C$ may change from line to line, and its dependence on the relevant parameters is indicated when needed. Finally, we set
\begin{equation}\label{AB}
\sigma=\frac{N(p-1)}{2}\in(0,p+1),
\end{equation}
so that $\sigma<2$, $\sigma=2$ and $\sigma>2$ correspond to the mass-subcritical, the mass-critical and the intercritical cases, respectively. Note that $p+1-\sigma=\frac{N+2-(N-2)p}{2}>0$ in the whole energy-subcritical range.

\smallskip
\noindent\textit{The effective coupling.} Throughout the paper, the damping coefficient is either a general function $\da\in C([0,\infty);[0,\infty))$, in which case $A$ is defined by \eqref{A}, or the polynomial damping $\da(t)=a(1+t)^{-\theta}$ of \eqref{Main-New}. In both cases, we set
\begin{equation}\label{def-h}
h(t)=\e^{-(p-1)A(t)},\qquad t\ge0 .
\end{equation}
Let $u$ be a solution of \eqref{general}, and set
\begin{equation}\label{change}
v(t,x)=\e^{A(t)}u(t,x).
\end{equation}
Since $i\partial_tu+\Delta u+i\da u=\e^{-A}(i\partial_tv+\Delta v)$ and $|u|^{p-1}u=\e^{-pA}|v|^{p-1}v$, the function $u$ solves \eqref{general} if and only if $v$ solves
\begin{equation}\label{V}
\begin{cases}
i\partial_t v+\Delta v+h(t)|v|^{p-1}v=0, & (t,x)\in(0,\infty)\times\R^N,\\
v(0)=u_0,
\end{cases}
\end{equation}
with $h$ given by \eqref{def-h}. The linear part of \eqref{general} generates the damped propagator $U_\da(t)=\e^{-A(t)}\e^{it\Delta}$. It is therefore natural to compare $u(t)$ with $U_\da(t)u_+$.

\begin{definition}\label{Def-Scat}
A global solution $u\in C([0,\infty);\Ha)$ of \eqref{general} \emph{scatters} (forward in time, in $\Ha$) if there exists $u_+\in\Ha$ such that
\begin{equation}\label{Scatt-def}
\lim_{t\to\infty}\big\|\e^{A(t)}\e^{-it\Delta}u(t)-u_+\big\|_{H^1}=0 .
\end{equation}
Equivalently, $\|\e^{-it\Delta}v(t)-u_+\|_{H^1}\to0$, or $\|u(t)-U_\da(t)u_+\|_{H^1}=\e^{-A(t)}o(1)$ as $t\to\infty$. We say that $u$ scatters in $L^2$ if \eqref{Scatt-def} holds with $u_+\in\Lt$ and with the $L^2$ norm.
\end{definition}

\begin{remark}\rm \label{rem:def}
(i) Since $\e^{-it\Delta}$ is unitary on $\Ha$, a solution that scatters satisfies $\|v(t)\|_{H^1}\to\|u_+\|_{H^1}$. In particular, we have $\sup_{t\ge0}\e^{A(t)}\|u(t)\|_{H^1}<\infty$.

(ii) When $A(t)\to\infty$, one always has $\|\e^{-it\Delta}u(t)\|_{L^2}=\|u(t)\|_{L^2}=m\e^{-A(t)}\to0$. Similarly, $\|\e^{-it\Delta}u(t)\|_{H^1}\to0$ for every solution with $\sup_t\e^{A(t)}\|u(t)\|_{H^1}<\infty$. Such statements carry no information on the asymptotic profile, and the weight $\e^{A(t)}$ in \eqref{Scatt-def} is essential. When $\theta>1$, we have $\e^{A(t)}\to\e^{A_\infty}$, where $A_\infty=\frac{a}{\theta-1}$, so that \eqref{Scatt-def} is equivalent to the usual scattering $\|\e^{-it\Delta}u(t)-\e^{-A_\infty}u_+\|_{H^1}\to0$.
\end{remark}

\noindent\textit{Admissible pairs and $p$-admissible exponents.} 

\begin{definition}\label{Admis}
A pair $(q,r)$ is \emph{admissible} if $\frac2q=N\big(\frac12-\frac1r\big)$ and $2\le r\le\frac{2N}{N-2}$ if $N\ge3$, $2\le r<\infty$ if $N=2$, $2\le r\le\infty$ if $N=1$. The admissible pair $(2,\frac{2N}{N-2})$, $N\ge3$, is called the endpoint.
\end{definition}

The scattering criterion below involves an auxiliary Lebesgue exponent $s$, which measures the norm in which the solution is a priori bounded, and an admissible pair $(q_s,r_s)$ tailored to $s$. Let $J_N=[2,\infty]$ if $N=1$, $J_N=[2,\infty)$ if $N=2$ and $J_N=[2,\frac{2N}{N-2}]$ if $N\ge3$, so that $\Ha\hookrightarrow L^s(\R^N)$ for $s\in J_N$.

\begin{definition}\label{def:padm}
Let $1<p<1+\frac4{(N-2)_+}$. An exponent $s\in J_N$ is called \emph{$p$-admissible} if the pair $(q_s,r_s)$ defined by
\begin{equation}\label{qsrs}
\frac1{r_s}=\frac12-\frac{p-1}{2s},\qquad\frac2{q_s}=\frac{N(p-1)}{2s}
\end{equation}
is admissible in the sense of Definition \ref{Admis}. We then define $k_s\in(1,\infty]$ by
\begin{equation}\label{ks}
\frac1{k_s}=1-\frac2{q_s}=1-\frac{N(p-1)}{2s}.
\end{equation}
\end{definition}

The following quantity measures the best possible time integrability exponent (see Lemma \ref{lem:ell} below):
\begin{equation}\label{ell}
\ell_N(p)=\begin{cases}\max\Big\{0,\,1-\dfrac{N(p-1)}4\Big\}, & N\ge2,\\[8pt]
\max\Big\{\dfrac12,\,1-\dfrac{p-1}4\Big\}, & N=1.\end{cases}
\end{equation}
In view of Lemmas \ref{lem:damping} and \ref{lem:ell} below, the following condition will be the key assumption for scattering in the polynomial case.

\begin{definition}\label{def:D}
We say that the damping $\da(t)=a(1+t)^{-\theta}$ satisfies condition $(\mathrm D_{N,p})$ if $a>0$ and either $0\le\theta<1$, or $\theta=1$ and $a(p-1)>\ell_N(p)$.
\end{definition}

\begin{remark}\rm \label{rem:Dexplicit}
For $\theta=1$, condition $(\mathrm D_{N,p})$ reads as follows. If $N\ge2$, it reduces to $p>1+\frac{4}{N+4a}$ in the mass-subcritical case, and to no condition other than $a>0$ if $p\ge1+\frac4N$. If $N=1$, it reads $p>1+\frac4{1+4a}$ and $a(p-1)>\frac12$. In particular, it reduces to $a>\frac18$ in the mass-critical case $p=5$. Indeed, for $p<1+\frac4N$, the inequality $a(p-1)>1-\frac{N(p-1)}4$ is equivalent to $(p-1)(N+4a)>4$.
\end{remark}

Finally, $K=K(N,p)$ denotes the optimal constant in the Gagliardo--Nirenberg inequality \eqref{Nirenberg}, and $Q$ denotes the ground state, that is, the unique positive radial solution of $-\Delta Q+Q=Q^{1+4/N}$ (see Lemma \ref{GNI}). We also write $p_{\rm s}(a)=1+\frac{4}{N+4a}$ and $p_{\rm ns}(a)=1+\frac{2}{N+2a}$.

\subsection{Main results}\label{S3}

\subsubsection*{Two abstract criteria}
We begin with two results for equation \eqref{V} with a general coupling $h$. They are the backbone of the paper, and they may be of independent interest. The first one is a scattering criterion.

\begin{theorem}[Scattering criterion]\label{thm:criterion}
Let $N\ge1$ and $1<p<1+\frac4{(N-2)_+}$. Let $h\in C([0,\infty))$ satisfy $0\le h\le1$, and let $v\in C([0,\infty);\Ha)$ be a global solution of \eqref{V} in the sense of Proposition \ref{prop:lwp}. Let $s$ be a $p$-admissible exponent, assume that
\begin{equation*}
L_s:=\sup_{t\ge0}\|v(t)\|_{L^s}<\infty,
\end{equation*}
and set $\Phi(t)=h(t)\|v(t)\|_{L^s}^{p-1}$. Assume that one of the following conditions holds:
\begin{enumerate}[(a)]
\item $q_s>2$ and $\Phi\in L^{k_s}(0,\infty)$;
\item $q_s=2$ (that is, $N\ge3$ and $s=\frac{N(p-1)}2$) and $\limsup_{t\to\infty}\Phi(t)<\delta_*:=(2pC_{\rm St})^{-1}$, where $C_{\rm St}$ is the constant of \eqref{St2}--\eqref{St3} for the endpoint pair.
\end{enumerate}
Then $v\in L^q(0,\infty;W^{1,r}(\R^N))$ for every admissible pair $(q,r)$, and there exists $v_+\in\Ha$ such that $\|\e^{-it\Delta}v(t)-v_+\|_{H^1}\to0$ as $t\to\infty$. In case (a), let $C_{\rm St}$ be the constant of Proposition \ref{prop:strichartz} for the pair $(q_s,r_s)$ and set $n=1+\big\lfloor\big(2pC_{\rm St}\|\Phi\|_{L^{k_s}(0,\infty)}\big)^{k_s}\big\rfloor$. Then
\begin{equation}\label{eq:Strbound}
\sup_{t\ge0}\|v(t)\|_{H^1}+\|v\|_{L^{q_s}(0,\infty;W^{1,r_s})}\le 2n\,(4C_{\rm St})^{n}\,\|u_0\|_{H^1},
\end{equation}
and, for all $t\ge0$,
\begin{equation}\label{rate}
\|\e^{-it\Delta}v(t)-v_+\|_{H^1}\le pC_{\rm St}\|\Phi\|_{L^{k_s}(t,\infty)}\|v\|_{L^{q_s}(t,\infty;W^{1,r_s})}.
\end{equation}
\end{theorem}

When $s=2$, the quantity $L_2=\|u_0\|_{L^2}$ is conserved, so that the criterion only involves the mass. The second result is a non-scattering criterion of Strauss--Barab type.

\begin{theorem}[Non-existence of asymptotically free solutions]\label{thm:nonscat}
Let $N\ge1$ and $1<p\le2$ with $p<1+\frac4{(N-2)_+}$. Let $h\in C([0,\infty))$ be nonnegative, and let $v\in C([0,\infty);\Ha)$ be a global solution of \eqref{V}. Assume that
\begin{equation}\label{nonint}
\int_1^\infty h(t)\,t^{-\frac{N(p-1)}2}\dd t=\infty .
\end{equation}
If there exists $v_+\in\Lt$ such that $\|\e^{-it\Delta}v(t)-v_+\|_{L^2}\to0$ as $t\to\infty$, then $v\equiv0$.
\end{theorem}

\begin{remark}\rm 
(i) For $h\equiv1$, Theorem \ref{thm:nonscat} recovers the classical fact that non-trivial solutions of the NLS equation are not asymptotically free in $L^2$ when $1<p\le1+\frac2N$ and $N\ge2$, or when $1<p\le2$ and $N=1$ \cite{Strauss1974}. The one-dimensional range $2<p\le3$ is due to Barab \cite{Barab1984}; see also \cite{TsutsumiYajima1984} and \cite[Section 7.5]{Cazenave2003}. The proof given in Section \ref{S6} only uses the conservation of mass and the factorisation \eqref{MDFM}. In particular, no weighted assumption on the data is needed, in contrast with the pseudo-conformal approach of \cite[Lemma 6.2]{BamriTayachi2023}, which requires $v_0\in\Sigma$. The conservation of mass used there holds for every real coefficient $h$, by the computation in the proof of Proposition \ref{prop:identities}.

(ii) Theorems \ref{thm:criterion} and \ref{thm:nonscat} are complementary. The former requires $h$ to be small at infinity in an $L^{k_s}$ sense. The latter requires $h$ not to decay faster than $t^{-1+N(p-1)/2}$, in an averaged sense. The gap between the two conditions reflects the fact that Theorem \ref{thm:criterion} does not use any dispersive decay of $v$. Section \ref{S-gap} indicates that, in the scale-invariant case, this gap is not merely technical.
\end{remark}

\subsubsection*{The polynomial damping: mass-subcritical case}
We now state our results for \eqref{Main-New}. Recall that $\sigma=\frac{N(p-1)}2$, $m=\|u_0\|_{L^2}$, that $K$ is the constant in \eqref{Nirenberg}, and that condition $(\mathrm D_{N,p})$ is given in Definition \ref{def:D}.

\begin{theorem}[Mass-subcritical case]\label{thm:sub}
Let $N\ge1$, $1<p<1+\frac4N$, $a\ge0$, $\theta\ge0$ and $u_0\in\Ha$.
\begin{enumerate}[(i)]
\item The solution $u$ of \eqref{Main-New} is global and, for all $t\ge0$,
\begin{equation}\label{bound-sub}
\e^{2A(t)}\|\nabla u(t)\|_{L^2}^2\le\max\Big\{2\|\nabla u_0\|_{L^2}^2,\ \Big(\frac{4K}{p+1}m^{p+1-\sigma}\Big)^{\frac{2}{2-\sigma}}\Big\},\qquad \e^{A(t)}\|u(t)\|_{L^2}=m .
\end{equation}
\item If $(\mathrm D_{N,p})$ holds, then $u$ scatters in $\Ha$. Explicitly, $u$ scatters if $a>0$ and $0\le\theta<1$. For $\theta=1$, it scatters if $p>1+\frac{4}{N+4a}$ when $N\ge2$, and if $p>1+\frac4{1+4a}$ and $a(p-1)>\frac12$ when $N=1$. Moreover, if $\theta=1$, $N\ge2$ and $p>1+\frac4{N+4a}$, the scattering state $u_+$ satisfies
\begin{equation}\label{rate-sub}
\big\|\e^{A(t)}\e^{-it\Delta}u(t)-u_+\big\|_{H^1}\le C(1+t)^{-\left((p-1)(a+\frac N4)-1\right)},\qquad t\ge0 ,
\end{equation}
with $C$ depending only on $N,p,a$ and $\|u_0\|_{H^1}$. The same estimate holds for $N=1$ when $1+\frac4{1+4a}<p\le3$.
\item Assume $u_0\neq0$ and $p\le2$. If either $\theta=1$ and $p\le1+\frac{2}{N+2a}$, or $\theta>1$ and $p\le1+\frac2N$, then $u$ does not scatter, even in $L^2$. In other words, there is no $u_+\in\Lt$ such that $\|\e^{A(t)}\e^{-it\Delta}u(t)-u_+\|_{L^2}\to0$.
\end{enumerate}
\end{theorem}

\begin{remark}\rm \label{rem:gap}
(i) Part (i) holds for every nonnegative continuous damping $\da(t)$ (see Proposition \ref{prop:apriori-sub}). It complements the global existence result of \cite[Proposition B.1]{HM-CPAA2023}, since it provides a bound on $\e^{A(t)}\|u(t)\|_{H^1}$ which is uniform in time and independent of the damping. For $\theta=0$, parts (i) and (ii) recover \cite[Theorem 2.2]{HM-MJM2025} in the mass-subcritical range, and they are consistent with the constant damping results of \cite{Inui2019}.

(ii) Let $\theta=1$ and $N\ge2$. Scattering holds for all data when $p>p_{\rm s}(a)$, while it fails for all non-trivial data when $p\le p_{\rm ns}(a)$. Note that $p_{\rm ns}(a)\le2$ when $N\ge2$. Since $\frac4{N+4a}>\frac2{N+2a}$, there is a gap. In this gap, Bamri and Tayachi \cite[Theorem 1.9]{BamriTayachi2023} constructed, for every $p-1>\frac{2}{N+2a}$, solutions with special oscillating data of the form $\e^{ib|x|^2/4}w_0$, where $b$ is sufficiently large and $w_0$ belongs to a weighted Sobolev space and satisfies the non-vanishing condition $\inf_x\langle x\rangle^n|w_0(x)|>0$ for a suitable integer $n$. These solutions scatter in $\Sigma$, hence in $\Ha$; see also \cite[Remark 1.10(1)]{BamriTayachi2023}. Consequently, $p_{\rm ns}(a)$ is the sharp threshold for the \emph{existence} of scattering solutions. Section \ref{S-gap} explains why $p_{\rm s}(a)$ should be the sharp threshold for the scattering of \emph{all} solutions (Conjecture \ref{conj:adiabatic}).

(iii) For $\theta=1$, part (iii) includes the borderline $p=p_{\rm ns}(a)$. This value is excluded in \cite[Theorem 1.7]{BamriTayachi2023}, since it cannot be decided under the sole assumption $t\da(t)\to\gamma$: by \cite[Propositions 7.1 and 7.2]{BamriTayachi2023}, logarithmic perturbations of $\gamma/(1+t)$ may lead to both behaviours at the critical value. For the pure power damping $\gamma/(1+t)$, \cite[Proposition 7.2]{BamriTayachi2023} (with $\tilde\gamma=0$) already shows that solutions with non-zero data in $\Sigma$ do not scatter at the critical value, under the additional condition $\gamma\ge1$ if $N=1$. Part (iii) extends this non-scattering result to all non-zero data in $\Ha$, without any weighted assumption. Note that, for $N=1$, our condition $p\le2$ at $p=p_{\rm ns}(a)$ is also equivalent to $\gamma=2a\ge1$.

(iv) For $\theta>1$ and $1+\frac2N<p<1+\frac4N$, we do not claim any scattering result for large data. Indeed, the damping is integrable, and the equation is asymptotically the undamped focusing NLS, for which solitary waves do not scatter. It is natural to expect that small data in $\Sigma$ scatter when $p-1$ is larger than the Strauss exponent $\alpha_1(0)$ (see \eqref{eq:quadratic}), by the arguments of \cite{Strauss1981,Tsutsumi1985,CazenaveWeissler1992}. Note that \cite[Theorem 1.2]{BamriTayachi2023} requires $0<\delta<1$, and therefore does not apply to a coupling $h$ bounded from below by $h_\infty>0$; the case $\delta=0$ is only mentioned in \cite[Remark 1.3(2)]{BamriTayachi2023}, with a reference to \cite{CazenaveWeissler1990}. We do not pursue this here.
\end{remark}

\subsubsection*{The polynomial damping: mass-critical case}
Let $p=1+\frac4N$, and let $Q$ be the ground state of Lemma \ref{GNI}. For a solution $u$ of \eqref{Main-New}, we define the time at which the mass of $u$ becomes subcritical,
\begin{equation}\label{TQ}
T_Q=\inf\big\{t\ge0:\ \e^{-A(t)}\|u_0\|_{L^2}<\|Q\|_{L^2}\big\}\in[0,\infty],
\end{equation}
with the convention $\inf\emptyset=\infty$. Since $A$ is continuous and nondecreasing, and strictly increasing when $a>0$, the time $T_Q$ can be computed explicitly from \eqref{A-explicit}. We have $T_Q=0$ if $\|u_0\|_{L^2}<\|Q\|_{L^2}$. If $\|u_0\|_{L^2}\ge\|Q\|_{L^2}$ and $a>0$, we set $L_Q=\log\big(\|u_0\|_{L^2}/\|Q\|_{L^2}\big)\ge0$, and we obtain
\begin{equation}\label{eq:TQ-explicit}
T_Q=\begin{cases}
\Big(1+\dfrac{(1-\theta)L_Q}{a}\Big)^{\frac1{1-\theta}}-1, & \theta<1, \ \text{ or } \ \theta>1 \text{ and } L_Q<\dfrac{a}{\theta-1},\\[10pt]
\Big(\dfrac{\|u_0\|_{L^2}}{\|Q\|_{L^2}}\Big)^{1/a}-1, & \theta=1,\\[10pt]
\infty, & \theta>1 \text{ and } L_Q\ge\dfrac{a}{\theta-1}.
\end{cases}
\end{equation}
Indeed, in this case $T_Q$ is the unique solution of $A(T_Q)=L_Q$ when $L_Q<\sup_{t\ge0}A(t)$, and $T_Q=\infty$ otherwise. In particular, $T_Q<\infty$ for every solution when $\theta\le1$ and $a>0$. For $\theta>1$, we have $T_Q<\infty$ if and only if $\e^{-a/(\theta-1)}\|u_0\|_{L^2}<\|Q\|_{L^2}$. When $a=0$, one has $T_Q=0$ if $\|u_0\|_{L^2}<\|Q\|_{L^2}$, and $T_Q=\infty$ otherwise.

\begin{theorem}[Mass-critical case]\label{thm:crit}
Let $N\ge1$, $p=1+\frac4N$, $a\ge0$, $\theta\ge0$ and $u_0\in\Ha$, and let $u$ be the maximal solution of \eqref{Main-New}.
\begin{enumerate}[(i)]
\item If $\|u_0\|_{L^2}<\|Q\|_{L^2}$, then $u$ is global and, with $\eta_0=(\|u_0\|_{L^2}/\|Q\|_{L^2})^{4/N}$,
\begin{equation}\label{bound-crit}
\sup_{t\ge0}\e^{2A(t)}\|\nabla u(t)\|_{L^2}^2\le\frac{E(u_0)}{1-\eta_0}\le\frac{\|\nabla u_0\|_{L^2}^2}{1-\eta_0}.
\end{equation}
\item More generally, if $T_{\max}>T_Q$, then $u$ is global and $\sup_{t\ge0}\e^{A(t)}\|u(t)\|_{H^1}<\infty$. In particular, if $\theta\le1$ and $a>0$, every global solution satisfies $\sup_{t\ge0}\e^{A(t)}\|u(t)\|_{H^1}<\infty$, and a solution is global as soon as it exists on $[0,T_Q]$.
\item If $(\mathrm D_{N,p})$ holds, that is, $a>0$ and $\theta<1$, or $\theta=1$ with $a>0$ if $N\ge2$ and $a>\frac18$ if $N=1$, then every global solution of \eqref{Main-New} scatters in $\Ha$.
\end{enumerate}
\end{theorem}

\begin{remark}\rm 
(i) Part (i) holds for every nonnegative continuous damping. It shows that the classical sufficient condition $\|u_0\|_{L^2}<\|Q\|_{L^2}$ for global existence remains valid for every time-dependent damping, whatever its decay. For $\theta=0$, part (iii) is \cite[Theorem 2.2]{HM-MJM2025} for $p=1+4/N$, which covers the range $1<p\le1+\frac4N$. Our proof is different, and it gives the more precise information (ii).

(ii) By (ii) and (iii), for $\theta\le1$, the only obstruction to global existence and scattering in the mass-critical case is a blow-up occurring before the explicit time $T_Q$ given by \eqref{eq:TQ-explicit}. Finite time blow-up does occur for some data when the damping is weak; see \cite{TsutsumiM1984,Darwich2012,HM-CPAA2023} for constant or general time-dependent dampings. Global existence under the assumption that $\Lsub$ is sufficiently large was obtained in \cite[Theorem 1.6]{HM-CPAA2023}.

(iii) For $\theta>1$ and $\e^{-a/(\theta-1)}\|u_0\|_{L^2}<\|Q\|_{L^2}$, part (ii) gives global existence with uniform $H^1$ bounds. Scattering in this case would require a time-dependent version of Dodson's theorem \cite{Dodson2015}, and we leave this question open.
\end{remark}

\subsubsection*{The polynomial damping: intercritical case}
Let $1+\frac4N<p<1+\frac{4}{(N-2)_+}$, so that $\sigma>2$. We define the scale-invariant quantity
\begin{equation}\label{Qfun}
\mathcal Q(f)=\|\nabla f\|_{L^2}^{\sigma-2}\|f\|_{L^2}^{p+1-\sigma},\qquad f\in\Ha,
\end{equation}
and the constants
\begin{equation}\label{eps0}
\varepsilon_0=\frac{p+1}{K\,2^{(\sigma+2)/2}},\qquad \kappa=\frac{p+1-\sigma}{\sigma-2}=\frac{N+2-(N-2)p}{N(p-1)-4}.
\end{equation}
Note that $\mathcal Q(u_0)^{1/(\sigma-2)}=\|\nabla u_0\|_{L^2}\|u_0\|_{L^2}^{\kappa}$ is the quantity appearing in \cite[(2.5)]{HM-MJM2025}, and that $\mathcal Q(f)^{1/(p-1)}=\|f\|_{L^2}^{1-s_c}\|\nabla f\|_{L^2}^{s_c}$, with $s_c=\frac N2-\frac2{p-1}$.

\begin{theorem}[Intercritical case, small data]\label{thm:inter-small}
Let $N\ge1$, $1+\frac4N<p<1+\frac{4}{(N-2)_+}$, $a\ge0$ and $\theta\ge0$, and let $u_0\in\Ha$ satisfy $\mathcal Q(u_0)<\varepsilon_0$. Then the solution $u$ of \eqref{Main-New} is global and
\begin{equation}\label{Decay-u}
\|\nabla u(t)\|_{L^2}\le\sqrt2\,\|\nabla u_0\|_{L^2}\,\e^{-A(t)},\qquad\|u(t)\|_{L^2}=\|u_0\|_{L^2}\e^{-A(t)},\qquad t\ge0.
\end{equation}
Moreover, $u$ scatters in $\Ha$ in each of the following cases:
\begin{enumerate}[(a)]
\item $(\mathrm D_{N,p})$ holds, that is, $a>0$ and $\theta<1$, or $\theta=1$ with $a>0$ if $N\ge2$ and $a(p-1)>\frac12$ if $N=1$;
\item $N\ge3$, $a\ge0$ and $\theta\ge0$ are arbitrary, and $\mathcal Q(u_0)<\varepsilon_1$, where $\varepsilon_1\in(0,\varepsilon_0]$ depends only on $N$ and $p$.
\end{enumerate}
\end{theorem}

\begin{theorem}[Intercritical case: scattering or grow-up]\label{thm:dichotomy}
Let $N\ge1$ and $1+\frac4N<p<1+\frac{4}{(N-2)_+}$. Assume $(\mathrm D_{N,p})$, and let $u$ be a global solution of \eqref{Main-New} with $u_0\neq0$. Then exactly one of the following alternatives holds.
\begin{enumerate}[(a)]
\item $u$ scatters in $\Ha$, and $\sup_{t\ge0}\e^{A(t)}\|u(t)\|_{H^1}<\infty$.
\item $\mathcal Q(u(t))\ge\varepsilon_0$ for every $t\ge0$. In this case,
\begin{equation}\label{growup}
\|\nabla u(t)\|_{L^2}\ge\Big(\frac{\varepsilon_0}{m^{p+1-\sigma}}\Big)^{\frac1{\sigma-2}}\e^{\kappa A(t)}\xrightarrow[t\to\infty]{}\infty ,
\end{equation}
and $\limsup_{t\to\infty}\|u(t)\|_{L^{p+1}}=\infty$.
\end{enumerate}
\end{theorem}

\begin{corollary}\label{cor:Lp+1}
Under the assumptions of Theorem \ref{thm:dichotomy}, $u$ scatters in $\Ha$ if one of the following conditions holds:
\begin{equation}\label{hyp-sup}
\sup_{t\ge0}\|u(t)\|_{L^{p+1}}<\infty;\qquad \liminf_{t\to\infty}\e^{-\kappa A(t)}\|\nabla u(t)\|_{L^2}=0;\qquad \sup_{t\ge0}\|u(t)\|_{H^1}<\infty .
\end{equation}
\end{corollary}

\begin{remark}\rm 
(i) For $\theta=0$, Theorem \ref{thm:inter-small} and Corollary \ref{cor:Lp+1} recover \cite[Theorems 2.3 and 2.4]{HM-MJM2025}. For $\theta\in(0,1]$, they are new, since $\Lsub=0$. As in \cite{HM-MJM2025}, the smallness condition only involves the scale-invariant quantity $\mathcal Q(u_0)$. It allows, for instance, $\|u_0\|_{L^2}\gg1$ with $\|\nabla u_0\|_{L^2}\ll1$.

(ii) Theorem \ref{thm:dichotomy} says that, when $\theta\le1$, a global solution which does not scatter must have a gradient that grows at least like $\e^{\kappa A(t)}$, up to constants. This means a growth like $\exp\big(\frac{\kappa a}{1-\theta}(1+t)^{1-\theta}\big)$ for $\theta<1$, and like $(1+t)^{\kappa a}$ for $\theta=1$. Whether such grow-up solutions exist is an open question. The mechanism behind the dichotomy is simple. The mass $\|u(t)\|_{L^2}=m\e^{-A(t)}$ tends to zero, so that a solution with bounded gradient eventually becomes small in the scale-invariant sense.

(iii) In case (b) of Theorem \ref{thm:inter-small} with $N\ge3$, scattering holds even for $\theta>1$ and for $a=0$. The proof relies on the endpoint Strichartz estimate and on the uniform bound \eqref{Decay-u}. For $N=1,2$ and $\theta>1$, we expect the same conclusion to hold by the small data theory in $\dot H^{s_c}$ (see \cite{CazenaveWeissler1990,HolmerRoudenko2008,Guevara2014}), adapted to the bounded coupling $h$. This is not proved here.
\end{remark}

\subsubsection*{More general dampings}
\begin{remark}\rm [More general dampings]\label{rem:general}
The proofs of Theorems \ref{thm:sub}--\ref{thm:dichotomy} only use the following properties of the damping: $\da\in C([0,\infty);[0,\infty))$, and $h=\e^{-(p-1)A}\in L^{k_s}(0,\infty)$ for some $p$-admissible $s$ with $q_s>2$. The latter property implies $A(t)\to\infty$, since $h$ is nonincreasing. The non-scattering statements only use the divergence \eqref{nonint}. These properties can be checked in terms of the quantity $t\,\da(t)$, which is the natural parameter in \cite{BamriTayachi2023}, as follows.
\end{remark}

\begin{corollary}[Dampings comparable to $1/t$]\label{cor:general}
Let $\da\in C([0,\infty);[0,\infty))$ and $1<p<1+\frac4{(N-2)_+}$.
\begin{enumerate}[(i)]
\item If $\liminf_{t\to\infty}t\,\da(t)>\frac{\ell_N(p)}{p-1}$, then the scattering conclusions of Theorems \ref{thm:sub}(ii) (except the rate \eqref{rate-sub}), \ref{thm:crit}(iii), \ref{thm:inter-small}(a) and \ref{thm:dichotomy}, as well as Corollary \ref{cor:Lp+1}, hold for \eqref{general}.
\item If $p\le2$ and there exist $b\ge0$ and $T\ge1$ such that $t\,\da(t)\le b$ for $t\ge T$ and $(p-1)(b+\frac N2)\le1$, then no non-trivial solution of \eqref{general} scatters in $L^2$.
\end{enumerate}
\end{corollary}

For instance, Corollary \ref{cor:general} applies to dampings $\da(t)=(1+t)^{-1}\chi(t)$ with $\chi$ continuous, nonnegative and $\chi(t)\to\infty$. It also applies to dampings which are regularly varying of index $-1$ in the sense of Karamata, that is, $\da(t)=t^{-1}L(t)$ for large $t$, with $L$ slowly varying and $\liminf_{t\to\infty}L(t)>\frac{\ell_N(p)}{p-1}$. Under the assumption $\Lsub>0$ of \cite{HM-MJM2025}, one has $h(t)\le\e^{-(p-1)\Lsub t}$, so that the scattering results of \cite{HM-MJM2025} are also recovered. In the notation of \cite{BamriTayachi2023}, where $\da_{\rm BT}=2\da$ and $t\,\da_{\rm BT}(t)\to\gamma$, condition (i) reads $\gamma>\frac{2\ell_N(p)}{p-1}$, and condition (ii) holds when $p-1\le\frac{2}{N+\gamma}$ and $t\,\da_{\rm BT}(t)\le\gamma$ for large $t$.

\subsubsection*{The gap in the scale-invariant case}
Our last main result describes the structure of the range $p_{\rm ns}(a)<p\le p_{\rm s}(a)$ left open by Theorem \ref{thm:sub} when $\theta=1$. It collects the rigorous statements of Section \ref{S-gap}. Let $\theta=1$, $a>0$ and $1<p<1+\frac4N$, and set
\begin{equation*}
\beta=\frac{2a(p-1)}{4-N(p-1)}>0,\qquad \tau(t)=\int_0^t(1+s)^{-2\beta}\dd s .
\end{equation*}
For $\omega>0$, let $R_\omega$ be the positive radial solution of $-\Delta R_\omega+\omega R_\omega=R_\omega^p$, and define the modulated ground state
\begin{equation*}
v_{\rm app}(t,x)=\e^{i\omega\tau(t)}\,(1+t)^{-\frac{N\beta}2}\exp\Big(\frac{i\beta|x|^2}{4(1+t)}\Big)R_\omega\Big(\frac{x}{(1+t)^{\beta}}\Big).
\end{equation*}

\begin{theorem}[Structure of the gap]\label{thm:gap}
Let $\theta=1$, $a>0$ and $1<p<1+\frac4N$.
\begin{enumerate}[(i)]
\item (Scaling.) The exponent $p_{\rm s}(a)$ is the unique exponent for which the natural scaling of the homogeneous problem $i\partial_tv+\Delta v+t^{-a(p-1)}|v|^{p-1}v=0$ leaves the $L^2$ norm invariant. Moreover, $p<p_{\rm s}(a)$ if and only if $\beta<\frac12$.
\item (Exact reduction.) A pseudo-conformal type change of variables, with spatial scale $(1+t)^{\beta}$ and time $\tau(t)$, transforms \eqref{Main-New} into the undamped focusing NLS equation \eqref{eq:repulsive} with a repulsive harmonic potential $q(\tau)|y|^2$. The coefficient $q$ is integrable and decays like $\tau^{-2}$ if $p<p_{\rm s}(a)$, and $q\equiv\frac1{16}$ if $p=p_{\rm s}(a)$.
\item (Concentration versus spreading.) If $\e^{A(t)}\|u(t)\|_{L^2(B(x(t),\rho(t)))}\ge c>0$ along a sequence $t\to\infty$, for some $\rho(t)=o(t)$, then $u$ does not scatter in $L^2$. Conversely, if $(p-1)(N+2a)\le4$ and $u_0\in\Sigma$, then $\sup_{x_0}\e^{A(t)}\|u(t)\|_{L^2(B(x_0,\rho))}\le C\big(\rho(1+t)^{-\beta}\big)^{\delta}$ for every $\rho>0$ and $t\ge1$, with $\delta>0$ and $C$ independent of $t$ and $\rho$.
\item (Approximate solitary waves.) For every $\omega>0$, $\|v_{\rm app}(t)\|_{L^2}=\|R_\omega\|_{L^2}$, and $v_{\rm app}$ solves \eqref{V} up to an error $\mathfrak r(t)$. If $p\neq1+\frac4{N+2a}$, then $\int_0^\infty\|\mathfrak r(t)\|_{L^2}\dd t<\infty$ if and only if $p<p_{\rm s}(a)$. If $p<1+\frac4{N+2a}$, then $v_{\rm app}$ is not asymptotically free.
\end{enumerate}
\end{theorem}

Theorem \ref{thm:gap} is a combination of Lemma \ref{lem:scaling}, Propositions \ref{prop:reduction}, \ref{prop:concentration}, \ref{prop:pseudoconf} and \ref{prop:approx}, where the precise statements and the explicit constants are given. It leads to Conjecture \ref{conj:adiabatic} below. This conjecture is supported by numerical simulations, and it asserts that non-scattering solutions exist in the open gap $(p_{\rm ns}(a),p_{\rm s}(a))$, so that $p_{\rm s}(a)$ would be sharp.

\subsection{Main ideas and comparison with the literature}
The guiding observation of this paper is that, for \eqref{Main-New}, the effective coupling $h(t)=\e^{-(p-1)A(t)}$ behaves as follows:
\begin{equation*}
h(t)=\begin{cases}
\exp\Big(-\dfrac{(p-1)a}{1-\theta}\big((1+t)^{1-\theta}-1\big)\Big), & 0\le\theta<1 \quad(\text{faster than any power}),\\[8pt]
(1+t)^{-a(p-1)}, & \theta=1\quad(\text{polynomial decay}),\\[6pt]
\searrow\ \e^{-(p-1)a/(\theta-1)}>0, & \theta>1\quad(\text{no decay}).
\end{cases}
\end{equation*}
The second observation is the elementary identity $h'=-(p-1)\da\,h$. It yields the exact relation
\begin{equation*}
h(t)+(p-1)\int_\tau^t \da(s)h(s)\dd s=h(\tau),\qquad 0\le\tau\le t,
\end{equation*}
which allows us to control the non-conserved energy of $v$ by the coupling at the initial time only. This simple fact replaces the bootstrap and Gr\"onwall-type arguments of \cite{HM-MJM2025}, and it makes all our a priori bounds independent of the damping. Once uniform $H^1$ bounds for $v$ are available, we derive scattering from a general criterion (Theorem \ref{thm:criterion}). This criterion only requires $h$ to belong to a suitable Lebesgue space $L^k(0,\infty)$, the exponent $k$ being dictated by the Strichartz estimates. We derive the non-scattering results from a direct version of the classical argument of Strauss and Barab, adapted to the coupling $h(t)$ (Theorem \ref{thm:nonscat}). This argument does not require weighted data.

\begin{figure}[ht]
\centering
\begin{tikzpicture}
\begin{axis}[width=12cm,height=7.2cm,xmin=0,xmax=2,ymin=1,ymax=2.45,
 xlabel={$\theta$},ylabel={$p$},
 xtick={0,0.5,1,1.5,2},ytick={1,1.4,1.5714,2.3333},
 yticklabels={$1$,$p_{\rm ns}(1)$,$p_{\rm s}(1)$,$1+\frac4N$},
 tick label style={font=\footnotesize},label style={font=\small}]
\addplot[fill=green!25,draw=none] coordinates {(0,1) (1,1) (1,2.3333) (0,2.3333)} \closedcycle;
\addplot[fill=red!20,draw=none] coordinates {(1,1) (2,1) (2,1.6667) (1,1.6667)} \closedcycle;
\addplot[pattern=north east lines,pattern color=black!35,draw=none] coordinates {(1,1.6667) (2,1.6667) (2,2.3333) (1,2.3333)} \closedcycle;
\addplot[line width=3pt,red!80!black] coordinates {(1,1) (1,1.4)};
\addplot[line width=3pt,orange!90!black] coordinates {(1,1.4) (1,1.5714)};
\addplot[line width=3pt,blue!70!black] coordinates {(1,1.5714) (1,2.3333)};
\addplot[thin,dashed,black!70] coordinates {(1,1.6667) (2,1.6667)};
\node[font=\footnotesize,anchor=south east] at (axis cs:2,1.6667) {$p=1+\frac2N$};
\addplot[thin,dotted,black!60] coordinates {(0,2.3333) (2,2.3333)};
\node[font=\footnotesize,align=center] at (axis cs:0.5,1.9) {all solutions\\ global and\\ scattering};
\node[font=\footnotesize,align=center] at (axis cs:1.5,1.3) {global, no\\ asymptotically\\ free solution};
\node[font=\footnotesize,align=center,fill=white,inner sep=1pt] at (axis cs:1.5,2.0) {global;\\ large data open};
\end{axis}
\end{tikzpicture}
\caption{The regimes of Theorem \ref{thm:sub} in the $(\theta,p)$-plane for $N=3$ and $a=1$, in the mass-subcritical range $1<p<1+\frac4N$. On the vertical segment $\theta=1$, all solutions scatter for $p>p_{\rm s}(1)=1+\frac47$ (blue part), and no non-trivial solution scatters for $p\le p_{\rm ns}(1)=1+\frac25$ (red part). The intermediate (orange) part is the gap studied in Section \ref{S-gap}. For $\theta>1$, no non-trivial solution scatters for $p\le1+\frac2N$, while the large data behaviour above $1+\frac2N$ is open.}
\label{fig:thetap}
\end{figure}
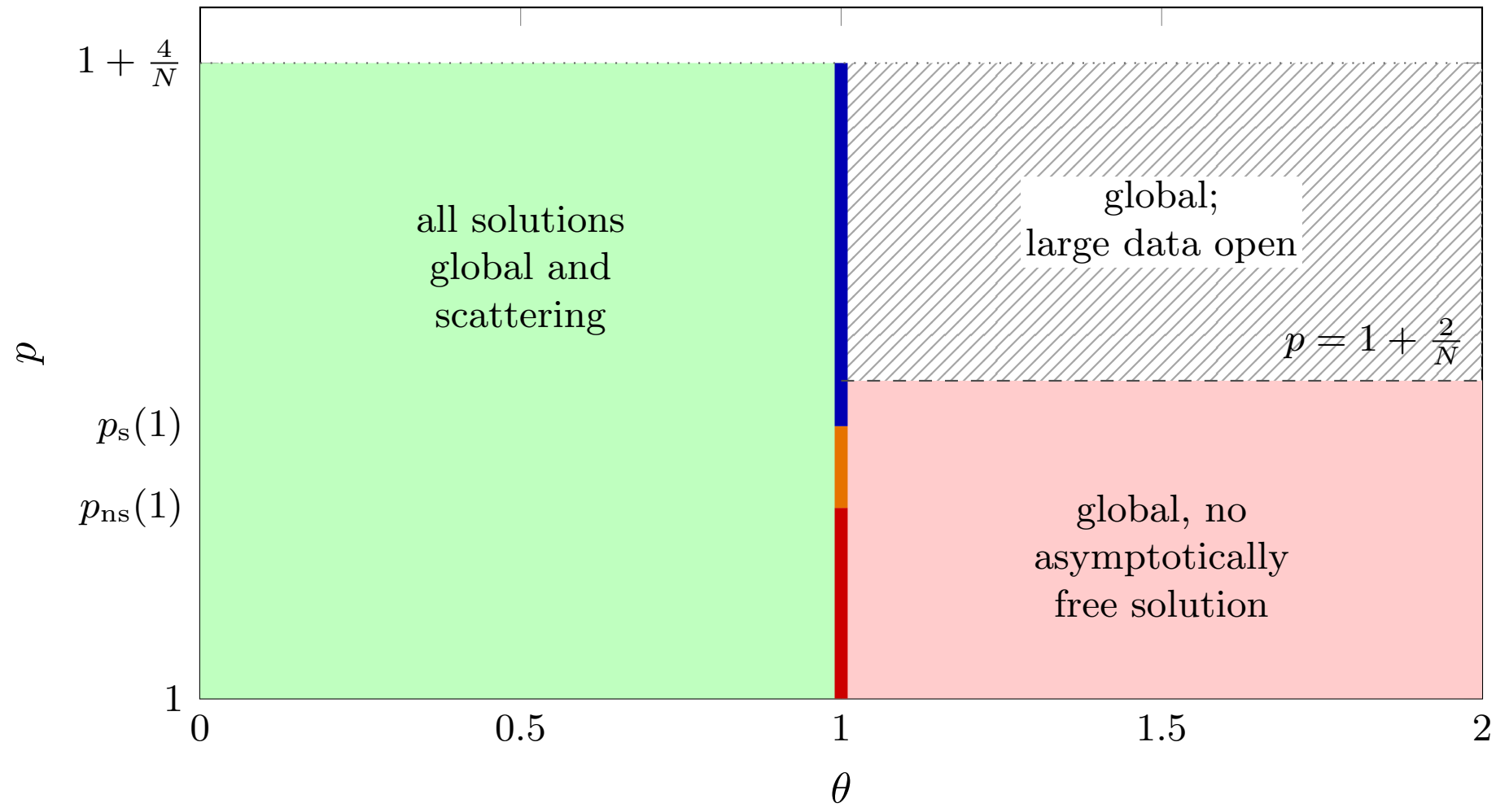

In the scale-invariant case $\theta=1$ and in the mass-subcritical regime, the picture is therefore the following. All solutions scatter above $p_{\rm s}(a)$, and none scatters below $p_{\rm ns}(a)$. In the intermediate range $(p_{\rm ns}(a),p_{\rm s}(a)]$, the results of \cite{BamriTayachi2023} provide scattering solutions for special oscillating data, while Section \ref{S-gap} indicates that non-scattering solutions should also exist there; see Figure \ref{fig:theta1} and Remark \ref{rem:gap}. Both thresholds decrease from the undamped values $1+4/N$ and $1+2/N$ to $1$ as $a\to\infty$. This quantifies the intuition that a stronger scale-invariant damping makes the nonlinearity more and more short-range. Table \ref{tab:summary} summarises our results, and Table \ref{tab:literature} compares them with the existing literature. Figure \ref{fig:thetap} displays the regimes in the $(\theta,p)$-plane.

\begin{figure}[ht]
\centering
\begin{tikzpicture}
\begin{axis}[width=11cm,height=7cm,xmin=0,xmax=5,ymin=1,ymax=2.45,
 xlabel={$a$},ylabel={$p$},
 legend style={at={(0.5,1.04)},anchor=south,font=\footnotesize,fill=white,legend columns=1},
 tick label style={font=\footnotesize},label style={font=\small}]
\addplot[domain=0:5,samples=200,thick,blue!70!black] {1+4/(3+4*x)};
\addlegendentry{$p_{\rm s}(a)=1+\frac{4}{3+4a}$ (all solutions scatter above)}
\addplot[domain=0:5,samples=200,thick,densely dotted,green!45!black] {1+(-(1+4*x)+sqrt((1+4*x)^2+16*(3+2*x)))/(2*(3+2*x))};
\addlegendentry{$p=1+\alpha_1(2a)$ (small data scatter above \cite{BamriTayachi2023})}
\addplot[domain=0:5,samples=200,thick,dashed,red!70!black] {1+2/(3+2*x)};
\addlegendentry{$p_{\rm ns}(a)=1+\frac{2}{3+2a}$ (no scattering on and below)}
\addplot[domain=0:5,thin,black!60,dotted] {1+4/3};
\addlegendentry{$p=1+4/N$}
\end{axis}
\end{tikzpicture}
\caption{The scale-invariant case $\theta=1$ in dimension $N=3$, mass-subcritical range. Above the solid curve, every $H^1$ solution scatters (Theorem \ref{thm:sub}). On and below the dashed curve, no non-trivial solution scatters (Theorem \ref{thm:sub}(iii)). Between the two curves, scattering is known for special oscillating data \cite{BamriTayachi2023}, and for small data in the Lorentz--Strichartz sense above the dotted curve \cite[Corollary 1.6]{BamriTayachi2023}, where $\alpha_1$ is defined by \eqref{eq:quadratic}. Conjecture \ref{conj:adiabatic} predicts the existence of non-scattering solutions everywhere below the solid curve.}
\label{fig:theta1}
\end{figure}
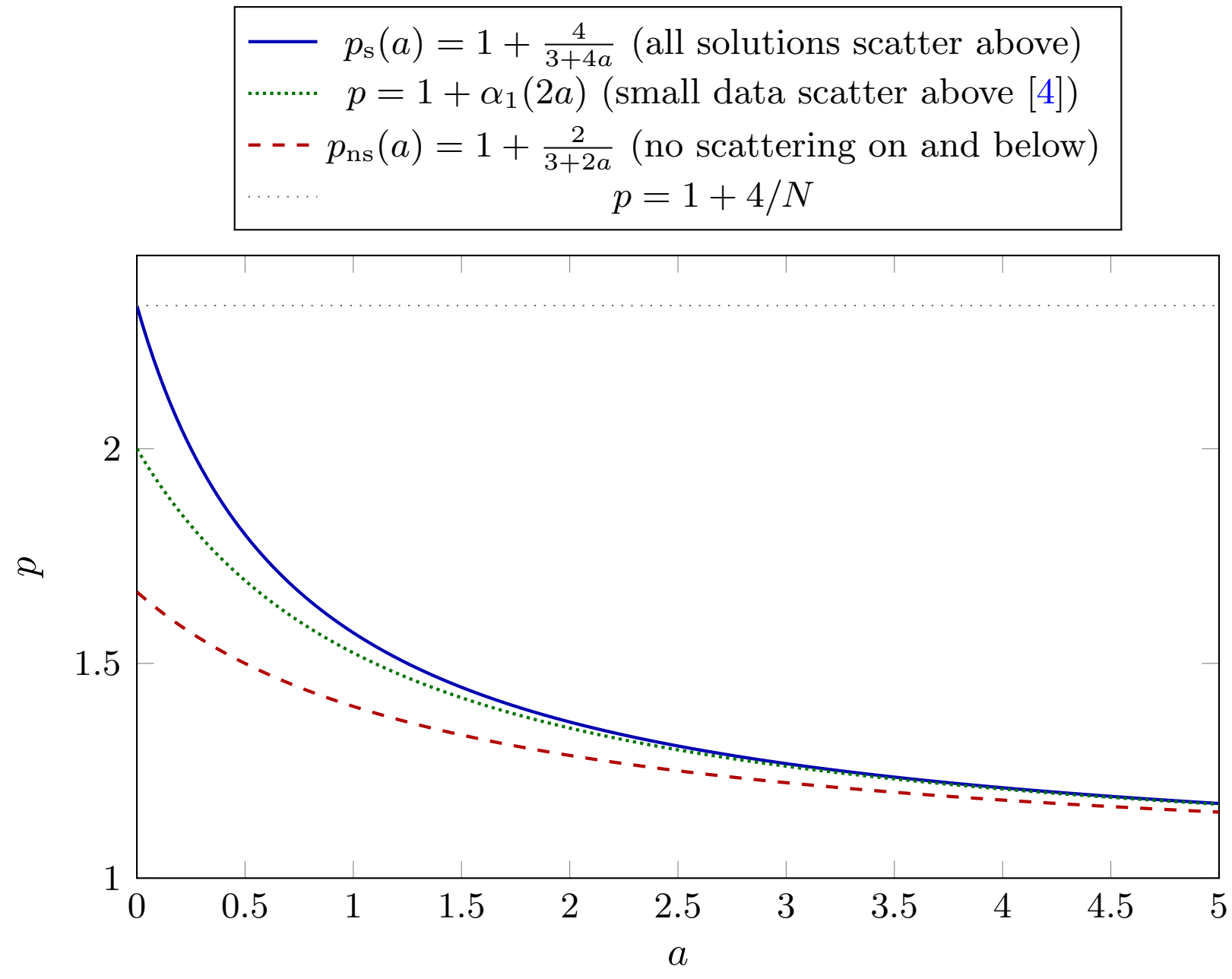

\begin{table}[ht]
\centering
\renewcommand{\arraystretch}{1.35}
{\small
\begin{tabular}{|p{3.2cm}|p{3.3cm}|p{3.6cm}|p{3.3cm}|}
\hline
 & $0\le\theta<1$ & $\theta=1$ & $\theta>1$\\
\hline
mass-subcritical\newline $1<p<1+\frac4N$ & all solutions global and scattering & all global; all scatter if $p>p_{\rm s}(a)$; none scatters if $p\le p_{\rm ns}(a)$ & all global; none scatters if $p\le 1+\frac2N$\\
\hline
mass-critical\newline $p=1+\frac4N$ & solutions existing beyond $T_Q$ are global and scatter & solutions existing beyond $T_Q$ are global and scatter & global if $\|u(\tau)\|_{L^2}<\|Q\|_{L^2}$ for some $\tau$\\
\hline
intercritical\newline $1+\frac4N<p$, $p<1+\frac{4}{(N-2)_+}$ & small data scatter; large data: scattering or grow-up & small data scatter; large data: scattering or grow-up & small data global; they scatter if $N\ge3$\\
\hline
\end{tabular}}
\caption{Summary of the results for \eqref{Main-New} with $a>0$. For $N=1$ and $\theta=1$, an additional lower bound on $a$ is needed, and the non-scattering statements require $p\le2$; see Section \ref{S3}.}
\label{tab:summary}
\end{table}

\begin{table}[ht]
\centering
\renewcommand{\arraystretch}{1.3}
{\footnotesize
\begin{tabular}{|p{2.6cm}|p{3.0cm}|p{4.1cm}|p{4.1cm}|}
\hline
Reference & Damping & Scattering & Non-scattering\\
\hline
Inui \cite{Inui2019}; \cite{HM-MJM2025} & constant, or $\Lsub>0$ & global solutions scatter in $H^1$ (under the conditions of \cite{HM-MJM2025}); small data (intercritical) & --\\
\hline
Bamri--Tayachi \cite{BamriTayachi2023} & $t\,\da_{\rm BT}(t)\to\gamma\ge0$, $L^1$-perturbations & small data (Lorentz--Strichartz), $N\ge2$; special oscillating data, $p-1>\frac2{N+\gamma}$ & data in $\Sigma$, $p-1<\frac2{N+\gamma}$, $p\le2$; at $p-1=\frac2{N+\gamma}$ for $\gamma/(1+t)$ (with $\gamma\ge1$ if $N=1$)\\
\hline
This paper, $\theta<1$ & $a(1+t)^{-\theta}$, $\Lsub=0$ & all $H^1$ data (mass-subcritical); beyond $T_Q$ (critical); small data and dichotomy (intercritical) & --\\
\hline
This paper, $\theta=1$ & $a(1+t)^{-1}$, $\gamma=2a$ & all $H^1$ data if $p>1+\frac4{N+4a}$ (mass-subcritical, $N\ge2$); mass-critical and intercritical as for $\theta<1$ & all $H^1$ data, $p\le1+\frac2{N+2a}$ (borderline included), $p\le2$; conjecturally some data for $p<1+\frac4{N+4a}$\\
\hline
This paper, $\theta>1$ & integrable & small data if $N\ge3$ (intercritical) & all $H^1$ data, $p\le1+\frac2N$, $p\le2$\\
\hline
\end{tabular}}
\caption{Comparison of the present results with the literature on the time-dependent damped NLS equation. The novelties are the treatment of $\Lsub=0$ for arbitrary $H^1$ data, the explicit scattering threshold $p_{\rm s}(a)$ and its interpretation, the non-scattering results without weighted assumptions (including the borderline exponent), and the scattering/grow-up dichotomy.}
\label{tab:literature}
\end{table}

\subsection{Organisation of the paper}
Section \ref{S2} collects the preliminaries: elementary properties of the damping and of the effective coupling, the local theory, the balance laws, the Gagliardo--Nirenberg and Strichartz estimates, and the Lebesgue exponents used in the scattering criterion. We prove the a priori bounds in Section \ref{S4}, the scattering criterion in Section \ref{S5}, and the non-scattering criterion in Section \ref{S6}. The proofs of Theorems \ref{thm:sub}--\ref{thm:dichotomy} and of Corollaries \ref{cor:Lp+1} and \ref{cor:general} are completed in Section \ref{S7}. Finally, Section \ref{S-gap} is devoted to the gap in the scale-invariant case, and contains the proof of Theorem \ref{thm:gap}.

\section{Preliminaries}\label{S2}

\subsection{The damping and the effective coupling}
Recall that $h=\e^{-(p-1)A}$ is defined by \eqref{def-h}. Since $A$ is nondecreasing and $A(0)=0$, the function $h$ is of class $C^1$ and nonincreasing. Moreover, it satisfies $0<h\le1$ and
\begin{equation}\label{h-prime}
h'(t)=-(p-1)\da(t)h(t),\qquad t\ge0.
\end{equation}
Integrating \eqref{h-prime} on $[\tau,t]$, we obtain the identity
\begin{equation}\label{eq:key}
h(t)+(p-1)\int_\tau^t \da(s)h(s)\dd s=h(\tau),\qquad 0\le\tau\le t,
\end{equation}
which plays a central role in Section \ref{S4}.

\begin{lemma}\label{lem:damping}
Let $\da(t)=a(1+t)^{-\theta}$ with $a>0$ and $\theta\ge0$, and let $A$ and $h$ be given by \eqref{A-explicit} and \eqref{def-h}.
\begin{enumerate}[(i)]
\item If $0\le\theta<1$, then $A(t)\to\infty$ as $t\to\infty$ and $h\in L^k(0,\infty)$ for every $k\in[1,\infty]$.
\item If $\theta=1$, then $A(t)=a\log(1+t)$ and $h(t)=(1+t)^{-a(p-1)}$. Moreover, $h\in L^k(0,\infty)$ for some $k\in[1,\infty)$ if and only if $k\,a(p-1)>1$, in which case
\begin{equation}\label{tail-h}
\|h\|_{L^k(t,\infty)}=\big(ka(p-1)-1\big)^{-1/k}\,(1+t)^{-a(p-1)+\frac1k},\qquad t\ge0.
\end{equation}
\item If $\theta>1$, then $A(t)\nearrow A_\infty=\frac{a}{\theta-1}$ and $h(t)\searrow h_\infty=\e^{-(p-1)a/(\theta-1)}>0$. In particular, $h\notin L^k(0,\infty)$ for every $k<\infty$.
\end{enumerate}
\end{lemma}

\begin{proof}
All assertions follow from \eqref{A-explicit}. In case (i), $A(t)=\frac{a}{1-\theta}\big((1+t)^{1-\theta}-1\big)$ grows like a positive power of $t$. It follows that $h(t)= \e^{(p-1)a/(1-\theta)}\exp\big(-\frac{(p-1)a}{1-\theta}(1+t)^{1-\theta}\big)$ decays faster than any power of $t$. Since $h\le1$, we deduce that $h\in L^k(0,\infty)$ for all $k\in[1,\infty]$. Case (ii) is an explicit computation, since $\int_t^\infty(1+s)^{-ka(p-1)}\dd s=\frac{(1+t)^{1-ka(p-1)}}{ka(p-1)-1}$ when $ka(p-1)>1$, while this integral diverges otherwise. Case (iii) follows from $A(t)=\frac{a}{\theta-1}\big(1-(1+t)^{1-\theta}\big)$.
\end{proof}

\begin{table}[ht]
\centering
\renewcommand{\arraystretch}{1.5}
{\small
\begin{tabular}{|c|>{\centering\arraybackslash}p{4.3cm}|>{\centering\arraybackslash}p{3.9cm}|>{\centering\arraybackslash}p{3.3cm}|}
\hline
 & $A(t)$ & $h(t)=\e^{-(p-1)A(t)}$ & integrability of $h$\\
\hline
$0\le\theta<1$ & $\frac{a}{1-\theta}\big((1+t)^{1-\theta}-1\big)\to\infty$ & exponential ($\theta=0$) or stretched exponential decay & $h\in L^k$ for all $k\in[1,\infty]$\\
\hline
$\theta=1$ & $a\log(1+t)\to\infty$ & $(1+t)^{-a(p-1)}$ & $h\in L^k\iff ka(p-1)>1$\\
\hline
$\theta>1$ & $\nearrow\frac{a}{\theta-1}<\infty$ & $\searrow h_\infty>0$ & $h\notin L^k$, $k<\infty$\\
\hline
\end{tabular}}
\caption{The three regimes of the polynomial damping $\da(t)=a(1+t)^{-\theta}$ (Lemma \ref{lem:damping}).}
\label{tab:regimes}
\end{table}

\begin{figure}[ht]
\centering
\begin{tikzpicture}
\begin{semilogyaxis}[width=11cm,height=6.2cm,xmin=0,xmax=60,ymin=1e-4,ymax=1.2,
 xlabel={$t$},ylabel={$h(t)$},
 legend style={at={(0.98,0.98)},anchor=north east,font=\footnotesize,fill=white},
 tick label style={font=\footnotesize},label style={font=\small}]
\addplot[domain=0:60,samples=300,thick,green!50!black] {exp(-2*((1+x)^0.5-1))};
\addlegendentry{$\theta=\frac12$}
\addplot[domain=0:60,samples=300,thick,orange!80!black,dashed] {(1+x)^(-1)};
\addlegendentry{$\theta=1$}
\addplot[domain=0:60,samples=300,thick,red!70!black,dotted] {exp(-2*(1-(1+x)^(-0.5)))};
\addlegendentry{$\theta=\frac32$}
\end{semilogyaxis}
\end{tikzpicture}
\caption{The effective coupling $h(t)$ for $a=1$ and $p=2$, in logarithmic scale: stretched exponential decay for $\theta<1$, polynomial decay for $\theta=1$, and convergence to the positive limit $\e^{-2}$ for $\theta>1$.}
\label{fig:h}
\end{figure}
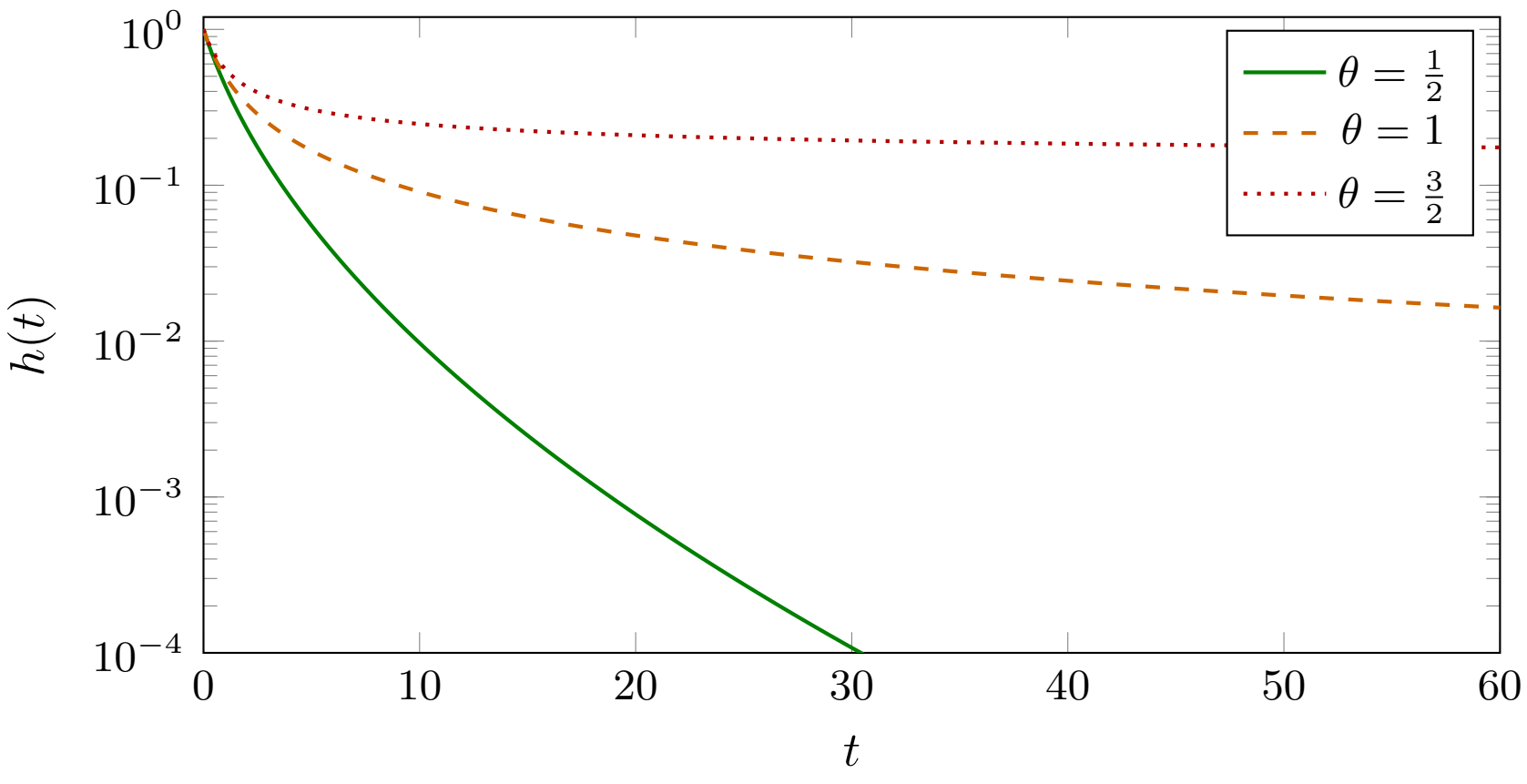

\subsection{Duhamel formula and local theory}
Let $v$ be given by \eqref{change}. In Duhamel form, equation \eqref{V} reads, for $0\le\tau\le t$ in the existence interval,
\begin{equation}\label{Integ-Eq}
v(t)=\e^{i(t-\tau)\Delta}v(\tau)+i\int_\tau^t\e^{i(t-s)\Delta}\big(h(s)|v(s)|^{p-1}v(s)\big)\dd s .
\end{equation}
Since $0<h\le1$ is continuous, the local theory of \eqref{V} in $\Ha$ is identical to that of the classical NLS equation. We recall it for the admissible pairs of Definition \ref{Admis}.

\begin{proposition}\label{prop:lwp}
Let $N\ge1$, $1<p<1+\frac{4}{(N-2)_+}$ and $\da\in C([0,\infty);[0,\infty))$. For every $u_0\in\Ha$, there exist $T_{\max}=T_{\max}(u_0)\in(0,\infty]$ and a unique maximal solution
\begin{equation*}
u\in C([0,T_{\max});\Ha)\cap L^q_{\rm loc}([0,T_{\max});W^{1,r}(\R^N))
\end{equation*}
of \eqref{general}, for every admissible pair $(q,r)$. The same holds for the function $v$ defined by \eqref{change}, which solves \eqref{V}. If $T_{\max}<\infty$, then $\|\nabla u(t)\|_{L^2}\to\infty$ as $t\uparrow T_{\max}$. Moreover, for every $\tau\in[0,T_{\max})$, the function $t\mapsto u(\tau+t)$ is the maximal solution of \eqref{general} with damping $\da(\tau+\cdot)$ and initial datum $u(\tau)$. The same conclusions hold for \eqref{V} with an arbitrary coefficient $h\in C([0,\infty))$ satisfying $0\le h\le1$.
\end{proposition}

\begin{proof}
The statement for \eqref{V} follows from the contraction argument of \cite[Theorem 4.4.6 and Proposition 4.4.1]{Cazenave2003} in Strichartz spaces. The coefficient $h$ only enters through the bound $\|h\|_{L^\infty}\le1$, so that the proof applies verbatim; see also \cite{Fanelli2009,GondaMachiharaOzawa2014}, where equations of the form \eqref{V} are treated for more general coefficients, and \cite[Remark 2.2]{HM-MJM2025}. In \cite{Cazenave2003}, the space-time integrability is stated for non-endpoint admissible pairs. For the endpoint pair ($N\ge3$), it follows on every compact interval $[0,T]\subset[0,T_{\max})$ by applying the endpoint Strichartz estimate \eqref{St2} below to \eqref{Integ-Eq}. Indeed, $h|v|^{p-1}v\in L^{\tilde q'}(0,T;W^{1,\tilde r'})$ for a suitable non-endpoint pair $(\tilde q,\tilde r)$, for instance $(\tilde q,\tilde r)=(\frac{4(p+1)}{N(p-1)},p+1)$. To see this, we use H\"older's inequality in space, $\|\,|v|^{p-1}\nabla v\|_{L^{(p+1)/p}}\le\|v\|_{L^{p+1}}^{p-1}\|\nabla v\|_{L^{p+1}}$, the bound $\sup_{[0,T]}\|v\|_{L^{p+1}}\le C\sup_{[0,T]}\|v\|_{H^1}$, and the inclusion $L^{\tilde q}(0,T)\subset L^{\tilde q'}(0,T)$, valid since $\tilde q>2$. The statement for $u$ follows from \eqref{change}, since $\e^{\pm A}$ is bounded on bounded time intervals. The last assertion is a consequence of uniqueness.
\end{proof}

\subsection{Balance laws}
For $f\in\Ha$, we use the notation of \cite{HM-MJM2025}:
\begin{equation}\label{Ener}
M(f)=\|f\|_{L^2}^2,\qquad E(f)=\|\nabla f\|_{L^2}^2-\frac{2}{p+1}\|f\|_{L^{p+1}}^{p+1},\qquad I(f)=\|\nabla f\|_{L^2}^2-\|f\|_{L^{p+1}}^{p+1}.
\end{equation}
For the solution $v$ of \eqref{V}, we introduce the Hamiltonian at time $t$,
\begin{equation}\label{Hamil-v}
\mathcal H(t)=\|\nabla v(t)\|_{L^2}^2-\frac{2}{p+1}h(t)\|v(t)\|_{L^{p+1}}^{p+1}=\e^{2A(t)}E(u(t)).
\end{equation}

\begin{proposition}\label{prop:identities}
Let $u\in C([0,T_{\max});\Ha)$ be the maximal solution of \eqref{general} and $v=\e^{A}u$. Then, for all $0\le\tau\le t<T_{\max}$,
\begin{align}
&\|v(t)\|_{L^2}=\|u_0\|_{L^2},\qquad\text{that is,}\qquad M(u(t))=\e^{-2A(t)}M(u_0),\label{M-Id}\\
&\mathcal H(t)=\mathcal H(\tau)+\frac{2(p-1)}{p+1}\int_\tau^t \da(s)h(s)\|v(s)\|_{L^{p+1}}^{p+1}\dd s,\label{Hamil-v-id}\\
&\frac{\dd}{\dd t}E(u(t))=-2\da(t)I(u(t)).\label{E-Id}
\end{align}
\end{proposition}

\begin{proof}
We first argue for sufficiently regular solutions. Writing \eqref{V} as $\partial_tv=i\Delta v+ih|v|^{p-1}v$, we obtain
\begin{equation*}
\frac{\dd}{\dd t}\|v\|_{L^2}^2=2\,{\rm Re}\,(\partial_tv,v)=2\,{\rm Re}\big(-i\|\nabla v\|_{L^2}^2+ih\|v\|_{L^{p+1}}^{p+1}\big)=0 .
\end{equation*}
Similarly, since $-\Delta v-h|v|^{p-1}v=i\partial_tv$, we get
\begin{equation*}
\begin{aligned}
\frac{\dd}{\dd t}\Big(\frac12\|\nabla v\|_{L^2}^2-\frac{h}{p+1}\|v\|_{L^{p+1}}^{p+1}\Big)
&={\rm Re}\big(-\Delta v-h|v|^{p-1}v,\partial_tv\big)-\frac{h'}{p+1}\|v\|_{L^{p+1}}^{p+1}\\
&={\rm Re}\big(i\partial_tv,\partial_tv\big)-\frac{h'}{p+1}\|v\|_{L^{p+1}}^{p+1},
\end{aligned}
\end{equation*}
and the first term on the right-hand side vanishes. Using \eqref{h-prime} and integrating in time, we obtain \eqref{Hamil-v-id}. Finally, \eqref{Hamil-v} gives $E(u(t))=\e^{-2A(t)}\mathcal H(t)$. Since $\e^{-2A}h\|v\|_{L^{p+1}}^{p+1}=\|u\|_{L^{p+1}}^{p+1}$, we deduce from \eqref{Hamil-v-id} and \eqref{h-prime} that
\begin{equation*}
\begin{aligned}
\frac{\dd}{\dd t}E(u)&=-2\da E(u)+\frac{2(p-1)}{p+1}\da\,\e^{-2A}h\|v\|_{L^{p+1}}^{p+1}\\
&=-2\da\Big(\|\nabla u\|_{L^2}^2-\frac{2}{p+1}\|u\|_{L^{p+1}}^{p+1}-\frac{p-1}{p+1}\|u\|_{L^{p+1}}^{p+1}\Big)=-2\da I(u).
\end{aligned}
\end{equation*}
The extension of these identities to $H^1$ solutions is standard. One may either regularise the data and use the continuous dependence in $H^1$, or apply directly the method of Ozawa \cite{Ozawa2006}, which avoids any regularisation and only uses the Duhamel formula. See also \cite{TsutsumiM1984}. The identities \eqref{M-Id} and \eqref{E-Id} are also stated in \cite[Proposition 3.1]{HM-MJM2025}, for sufficiently smooth solutions. In the regularisation approach, $H^2$ data give $H^2$ solutions on the $H^1$ existence interval by Kato's method \cite{Kato1987}, \cite[Section 5.3]{Cazenave2003}, which consists in estimating $\partial_tv$ in Strichartz norms. For \eqref{V}, differentiating the equation in time only produces the additional term $h'|v|^{p-1}v$, which is harmless because $h'$ is bounded on compact intervals. The computation above is then justified for such solutions. The identities pass to the limit by continuous dependence in $H^1$, the right-hand side of \eqref{Hamil-v-id} being continuous in $H^1$ by the Sobolev embedding. Finally, since the integrand in \eqref{Hamil-v-id} is continuous in time, $\mathcal H$ is of class $C^1$, and so is $E(u(t))=\e^{-2A(t)}\mathcal H(t)$; hence \eqref{E-Id} holds pointwise.
\end{proof}

\begin{lemma}\label{lem:liminf}
Let $u$ be a global solution of \eqref{general} with $\int_0^\infty \da(s)\dd s=\infty$. Then $\liminf_{t\to\infty}I(u(t))\le0$.
\end{lemma}

\begin{proof}
We argue by contradiction and assume that there exist $\delta>0$ and $T_0\ge0$ such that $I(u(t))\ge\delta$ for $t\ge T_0$. By \eqref{E-Id}, we get $E(u(t))\le E(u(T_0))-2\delta\int_{T_0}^t\da(s)\dd s$, and the right-hand side tends to $-\infty$. On the other hand, $E(f)\ge I(f)$ for every $f$, since $\frac{2}{p+1}<1$. We deduce that $E(u(t))\ge\delta>0$ for $t\ge T_0$, which is a contradiction. This argument is taken from \cite[Lemma 3.7]{HM-MJM2025}, where it is stated for $N\ge3$; the proof does not depend on the dimension.
\end{proof}

\subsection{Gagliardo--Nirenberg inequalities}
Let $1<p<1+\frac4{(N-2)_+}$, and recall that $\sigma=\frac{N(p-1)}2\in(0,p+1)$ is defined by \eqref{AB}.

\begin{lemma}\label{GNI}
There exists an optimal constant $K=K(N,p)>0$ such that
\begin{equation}\label{Nirenberg}
\|f\|_{L^{p+1}}^{p+1}\le K\|f\|_{L^2}^{p+1-\sigma}\|\nabla f\|_{L^2}^{\sigma},\qquad f\in\Ha .
\end{equation}
For $p=1+\frac4N$, one has $K=\frac{p+1}{2}\|Q\|_{L^2}^{-4/N}$, where $Q$ is the unique positive radial solution of $-\Delta Q+Q=Q^{1+4/N}$ \cite{Weinstein1983,Kwong1989}. In other words,
\begin{equation}\label{eq:sharpGN}
\frac{2}{p+1}\|f\|_{L^{p+1}}^{p+1}\le\Big(\frac{\|f\|_{L^2}}{\|Q\|_{L^2}}\Big)^{4/N}\|\nabla f\|_{L^2}^2,\qquad p=1+\frac4N .
\end{equation}
Moreover, in the intercritical case $1+\frac4N<p<1+\frac{4}{(N-2)_+}$, setting $s_*=\frac{N(p-1)}2\in(2,2^*)$ (note that $s_*=\sigma$) and $s_c=\frac N2-\frac2{p-1}\in(0,1)$, there exists $K_*=K_*(N,p)$ such that
\begin{equation}\label{eq:GNcrit}
\|f\|_{L^{s_*}}\le K_*\|f\|_{L^2}^{1-s_c}\|\nabla f\|_{L^2}^{s_c},\qquad f\in\Ha .
\end{equation}
\end{lemma}

\begin{proof}
Inequalities \eqref{Nirenberg} and \eqref{eq:GNcrit} are particular cases of the Gagliardo--Nirenberg inequality \cite{Cazenave2003,LinaresPonce2015}, and the exponents are dictated by scaling. Inequality \eqref{eq:GNcrit} also follows from the Sobolev embedding $\dot H^{s_c}\hookrightarrow L^{s_*}$ (note that $s_c<\frac N2$) and the interpolation inequality $\|f\|_{\dot H^{s_c}}\le\|f\|_{L^2}^{1-s_c}\|\nabla f\|_{L^2}^{s_c}$. The value of the optimal constant in the mass-critical case is due to Weinstein \cite{Weinstein1983}, and the uniqueness of $Q$ is due to Kwong \cite{Kwong1989}. For the reader's convenience, we note that the value of $K$ can be checked on $Q$ itself: the Pohozaev identity gives $\|\nabla Q\|_{L^2}^2=\frac{2}{p+1}\|Q\|_{L^{p+1}}^{p+1}$, so that equality holds in \eqref{eq:sharpGN} for $f=Q$.
\end{proof}

\subsection{Strichartz and dispersive estimates}
We shall use the following classical estimates; see \cite{Strichartz1977,GinibreVelo1985,Yajima1987,KeelTao1998} and \cite[Chapter 2]{Cazenave2003}.

\begin{proposition}\label{prop:strichartz}
Let $(q,r)$ and $(\tilde q,\tilde r)$ be admissible pairs. There exists a constant $C_{\rm St}=C_{\rm St}(N,q,r,\tilde q,\tilde r)\ge1$ such that, for every interval $J\ni\tau$, every $\varphi\in\Lt$ and every $F\in L^{\tilde q'}(J;L^{\tilde r'}(\R^N))$,
\begin{align}
&\|\e^{it\Delta}\varphi\|_{L^q(\R;L^r)}\le C_{\rm St}\|\varphi\|_{L^2},\label{St1}\\
&\Big\|\int_\tau^t\e^{i(t-s)\Delta}F(s)\dd s\Big\|_{L^q(J;L^r)\cap L^\infty(J;L^2)}\le C_{\rm St}\|F\|_{L^{\tilde q'}(J;L^{\tilde r'})},\label{St2}\\
&\Big\|\int_J\e^{-is\Delta}F(s)\dd s\Big\|_{L^2}\le C_{\rm St}\|F\|_{L^{\tilde q'}(J;L^{\tilde r'})}.\label{St3}
\end{align}
Here the norm of $L^q\cap L^\infty$ is the sum of the two norms. The endpoint pair is allowed when $N\ge3$ \cite{KeelTao1998}. Moreover, for $2\le r\le\infty$ and $t\neq0$,
\begin{equation}\label{dispersive}
\|\e^{it\Delta}\varphi\|_{L^r}\le(4\pi|t|)^{-N(\frac12-\frac1r)}\|\varphi\|_{L^{r'}} .
\end{equation}
\end{proposition}

Since $\nabla$ commutes with $\e^{it\Delta}$, the estimates \eqref{St1}--\eqref{St3} also hold when $L^r$, $L^{\tilde r'}$ and $L^2$ are replaced by $W^{1,r}$, $W^{1,\tilde r'}$ and $H^1$ (with the norms fixed in Subsection \ref{S1-setting}), up to replacing $C_{\rm St}$ by $2C_{\rm St}$. In the sequel, $C_{\rm St}$ denotes this larger constant, so that \eqref{St1}--\eqref{St3} may be used indifferently in $L^r$ or in $W^{1,r}$.

We also recall the factorisation of the Schr\"odinger group. For $t>0$, let $\Mo(t)$ be the multiplication by $\e^{i|x|^2/(4t)}$, and let $\Do(t)$ be the dilation $\Do(t)g(x)=(2it)^{-N/2}g\big(\frac{x}{2t}\big)$. Then, for $\varphi\in\Lt$ and $t>0$, we have
\begin{equation}\label{MDFM}
\e^{it\Delta}\varphi=\Mo(t)\Do(t)\Fou\big(\Mo(t)\varphi\big),\qquad \|\Do(t)g\|_{L^r}=(2t)^{-N(\frac12-\frac1r)}\|g\|_{L^r}.
\end{equation}
The first identity is a rewriting of the explicit kernel $(4\pi it)^{-N/2}\e^{i|x-y|^2/(4t)}$; see \cite[Section 7.2]{Cazenave2003}. In particular, $\Mo(t)$ and $\Do(t)$ are unitary on $\Lt$. (We use calligraphic letters to avoid any confusion with the mass functional $M$ in \eqref{Ener}.)

\subsection{Lebesgue exponents for the scattering criterion}
We now study the $p$-admissible exponents of Definition \ref{def:padm}. The definition is designed so that, by H\"older's inequality,
\begin{equation}\label{Holder-rs}
\frac1{r_s'}=\frac{p-1}{s}+\frac1{r_s},\qquad\frac1{q_s'}=\frac1{k_s}+\frac1{q_s}.
\end{equation}

\begin{lemma}\label{lem:ell}
Let $N\ge1$ and $1<p<1+\frac4{(N-2)_+}$. For every $\varepsilon>0$, there exists a $p$-admissible exponent $s$ with $q_s>2$ and $\frac1{k_s}<\ell_N(p)+\varepsilon$. If $\ell_N(p)>0$, one can take $\varepsilon=0$, that is, there is a $p$-admissible $s$ with $q_s>2$ and $\frac1{k_s}=\ell_N(p)$. Moreover, $\frac1{k_s}\ge\ell_N(p)$ for every $p$-admissible $s$.
\end{lemma}

\begin{proof}
Given $s\in J_N$, the relations \eqref{qsrs} automatically yield $\frac2{q_s}=N\big(\frac12-\frac1{r_s}\big)$ and $r_s\ge2$. Hence $s$ is $p$-admissible if and only if $r_s$ lies in the range of Definition \ref{Admis}, that is, $s\ge\frac{N(p-1)}2$ if $N\ge3$, $s>p-1$ if $N=2$, and $s\ge p-1$ if $N=1$. In all cases, $\frac1{k_s}=1-\frac{N(p-1)}{2s}$ is increasing in $s$. It follows that the infimum of $\frac1{k_s}$ is obtained for the smallest admissible $s$. Moreover, $q_s>2$ means $s>\frac{N(p-1)}2$.

Let $N\ge3$. If $p<1+\frac4N$, the choice $s=2$ is $p$-admissible, since $2>\frac{N(p-1)}2$. It gives $q_s>2$ and $\frac1{k_s}=1-\frac{N(p-1)}4=\ell_N(p)>0$. If $p\ge1+\frac4N$, then $\ell_N(p)=0$. The exponents $s\in\big(\frac{N(p-1)}{2},\frac{2N}{N-2}\big]$ are $p$-admissible with $q_s>2$, this interval being non-empty because $p<1+\frac4{N-2}$, and $\frac1{k_s}\to0$ as $s\downarrow\frac{N(p-1)}2$. Let $N=2$. If $p<3$, then $s=2$ is $p$-admissible with $q_s>2$ and $\frac1{k_s}=1-\frac{p-1}2=\ell_2(p)$. If $p\ge3$, every $s\in(p-1,\infty)$ is $p$-admissible with $q_s>2$, and $\frac1{k_s}\to0$ as $s\downarrow p-1$. Let finally $N=1$. If $p\le3$, then $s=2$ is $p$-admissible, $q_s=\frac8{p-1}\ge4$ and $\frac1{k_s}=1-\frac{p-1}4=\ell_1(p)$. If $p>3$, then $s=p-1$ is $p$-admissible with $r_s=\infty$, $q_s=4$ and $\frac1{k_s}=\frac12=\ell_1(p)$. The last assertion follows from the characterisation of $p$-admissibility given above and from the monotonicity of $s\mapsto\frac1{k_s}$, since $s\ge2$ in all cases and $s\ge p-1$ when $N=1$.
\end{proof}

\begin{table}[ht]
\centering
\renewcommand{\arraystretch}{1.5}
{\small
\begin{tabular}{|c|c|c|c|c|}
\hline
 & range of $p$ & exponent $s$ & pair $(q_s,r_s)$ & $\frac1{k_s}$\\
\hline
$N\ge2$ & $p<1+\frac4N$ & $s=2$ (mass only) & $\big(\frac{8}{N(p-1)},\frac{4}{3-p}\big)$ & $1-\frac{N(p-1)}4=\ell_N(p)$\\
\hline
$N\ge3$ & $p\ge1+\frac4N$ & $s\downarrow\frac{N(p-1)}2$ & $\to\big(2,\frac{2N}{N-2}\big)$ & $\to0$\\
\hline
$N=2$ & $p\ge3$ & $s\downarrow p-1$ & $q_s\downarrow2$, $r_s\uparrow\infty$ & $\to0$\\
\hline
$N=1$ & $p\le3$ & $s=2$ (mass only) & $\big(\frac8{p-1},\frac4{3-p}\big)$ & $1-\frac{p-1}4$\\
\hline
$N=1$ & $p>3$ & $s=p-1$ & $(4,\infty)$ & $\frac12$\\
\hline
\end{tabular}}
\caption{Optimal choices of the exponent $s$ in Lemma \ref{lem:ell}. Scattering for $\theta=1$ requires $h\in L^{k_s}(0,\infty)$, that is, $a(p-1)>\frac1{k_s}$.}
\label{tab:exponents}
\end{table}

\begin{lemma}\label{lem:Dimplies}
If $(\mathrm D_{N,p})$ holds, then $A(t)\to\infty$ as $t\to\infty$, and there exists a $p$-admissible exponent $s$ with $q_s>2$ such that $h\in L^{k_s}(0,\infty)$.
\end{lemma}

\begin{proof}
If $\theta<1$, any $p$-admissible $s$ with $q_s>2$ will do, by Lemma \ref{lem:damping}(i). If $\theta=1$, Lemma \ref{lem:ell} with $\varepsilon=a(p-1)-\ell_N(p)>0$ provides a $p$-admissible $s$ with $q_s>2$ and $\frac1{k_s}<a(p-1)$. Then $k_sa(p-1)>1$, and Lemma \ref{lem:damping}(ii) gives $h\in L^{k_s}(0,\infty)$. In both cases, $A(t)\to\infty$ by Lemma \ref{lem:damping}.
\end{proof}

\section{Uniform energy bounds}\label{S4}

In this section, $\da\in C([0,\infty);[0,\infty))$ is an arbitrary nonnegative damping, and $u\in C([0,T_{\max});\Ha)$ is the maximal solution of \eqref{general}. We set $v=\e^{A}u$, $h=\e^{-(p-1)A}$ and $m=\|u_0\|_{L^2}=\|v(t)\|_{L^2}$. For $0\le\tau\le t<T_{\max}$, we write
\begin{equation*}
f(t)=\|\nabla v(t)\|_{L^2}^2,\qquad F_\tau(t)=\sup_{\tau\le s\le t}f(s).
\end{equation*}
All the bounds below rest on the following consequence of \eqref{Hamil-v} and \eqref{Hamil-v-id}: for $0\le\tau\le s<T_{\max}$, we have
\begin{equation}\label{basic}
f(s)=\mathcal H(\tau)+\frac2{p+1}\Big(h(s)\|v(s)\|_{L^{p+1}}^{p+1}+(p-1)\int_\tau^s\da(\varsigma)h(\varsigma)\|v(\varsigma)\|_{L^{p+1}}^{p+1}\dd\varsigma\Big).
\end{equation}
Assume that $\|v(\varsigma)\|_{L^{p+1}}^{p+1}\le G$ for $\varsigma\in[\tau,s]$. Thanks to the identity \eqref{eq:key}, the bracket is then bounded by $\big(h(s)+(p-1)\int_\tau^s\da h\big)G=h(\tau)G$. This observation is exact, and it loses nothing with respect to the undamped case. It is the main point of the following three propositions.

\begin{proposition}[Mass-subcritical case]\label{prop:apriori-sub}
Let $1<p<1+\frac4N$. Then $T_{\max}=\infty$ and, for all $t\ge0$,
\begin{equation}\label{bound-sub2}
f(t)\le\max\Big\{2\|\nabla u_0\|_{L^2}^2,\ \Big(\frac{4K}{p+1}m^{p+1-\sigma}\Big)^{\frac{2}{2-\sigma}}\Big\}.
\end{equation}
\end{proposition}

\begin{proof}
Let $0\le s\le t<T_{\max}$. By \eqref{Nirenberg}, we have $\|v(\varsigma)\|_{L^{p+1}}^{p+1}\le Km^{p+1-\sigma}f(\varsigma)^{\sigma/2}\le Km^{p+1-\sigma}F_0(t)^{\sigma/2}$ for $\varsigma\in[0,t]$. We insert this bound in \eqref{basic} with $\tau=0$. Using \eqref{eq:key} with $h(0)=1$ and the bound $\mathcal H(0)=E(u_0)\le\|\nabla u_0\|_{L^2}^2$, we obtain $f(s)\le\|\nabla u_0\|_{L^2}^2+BF_0(t)^{\sigma/2}$, where $B=\frac{2K}{p+1}m^{p+1-\sigma}$ and $\frac\sigma2\in(0,1)$. Taking the supremum over $s\in[0,t]$, we get
\begin{equation*}
F_0(t)\le\|\nabla u_0\|_{L^2}^2+BF_0(t)^{\sigma/2} .
\end{equation*}
If $F_0(t)\ge(2B)^{2/(2-\sigma)}$, then $BF_0(t)^{\sigma/2}=BF_0(t)^{-(2-\sigma)/2}F_0(t)\le\frac12F_0(t)$, and the previous inequality gives $F_0(t)\le2\|\nabla u_0\|_{L^2}^2$. Hence, in all cases, $F_0(t)$ is bounded by the right-hand side of \eqref{bound-sub2}, which does not depend on $t$. Since $\|\nabla u(t)\|_{L^2}=\e^{-A(t)}f(t)^{1/2}\le f(t)^{1/2}$, the blow-up alternative of Proposition \ref{prop:lwp} gives $T_{\max}=\infty$. We conclude that \eqref{bound-sub2} holds for all $t\ge0$.
\end{proof}

\begin{proposition}[Mass-critical case]\label{prop:apriori-crit}
Let $p=1+\frac4N$, and assume that there exists $\tau\in[0,T_{\max})$ such that $\|u(\tau)\|_{L^2}<\|Q\|_{L^2}$. Set $\eta_\tau=(\|u(\tau)\|_{L^2}/\|Q\|_{L^2})^{4/N}\in[0,1)$. Then $T_{\max}=\infty$ and
\begin{equation}\label{bound-crit2}
\sup_{t\ge\tau}f(t)\le\frac{\mathcal H(\tau)}{1-\eta_\tau}\le\frac{f(\tau)}{1-\eta_\tau}.
\end{equation}
\end{proposition}

\begin{proof}
Since $\|u(\tau)\|_{L^2}=\e^{-A(\tau)}m$ and $h(\tau)=\e^{-\frac4NA(\tau)}$, we have $\eta_\tau=(m/\|Q\|_{L^2})^{4/N}h(\tau)$. By \eqref{eq:sharpGN} and the conservation of $\|v\|_{L^2}=m$, we get $\frac2{p+1}\|v(\varsigma)\|_{L^{p+1}}^{p+1}\le(m/\|Q\|_{L^2})^{4/N}f(\varsigma)$. Inserting this bound in \eqref{basic} and using \eqref{eq:key}, we obtain, for $\tau\le s\le t<T_{\max}$,
\begin{equation*}
f(s)\le\mathcal H(\tau)+\Big(\frac m{\|Q\|_{L^2}}\Big)^{4/N}\Big(h(s)+(p-1)\int_\tau^s\da h\Big)F_\tau(t)=\mathcal H(\tau)+\eta_\tau F_\tau(t).
\end{equation*}
Taking the supremum over $s\in[\tau,t]$, we infer that $F_\tau(t)\le\mathcal H(\tau)+\eta_\tau F_\tau(t)$, that is, $F_\tau(t)\le\mathcal H(\tau)/(1-\eta_\tau)$. Applying \eqref{eq:sharpGN} once more at time $\tau$, we get $(1-\eta_\tau)f(\tau)\le\mathcal H(\tau)\le f(\tau)$. In particular, $\mathcal H(\tau)\ge0$. Since the bound does not depend on $t$, the blow-up alternative gives $T_{\max}=\infty$.
\end{proof}

\begin{proposition}[Intercritical case]\label{prop:apriori-inter}
Let $1+\frac4N<p<1+\frac4{(N-2)_+}$, and assume that there exists $\tau\in[0,T_{\max})$ such that $\mathcal Q(u(\tau))<\varepsilon_0$, where $\mathcal Q$ and $\varepsilon_0$ are defined in \eqref{Qfun}--\eqref{eps0}. Then $T_{\max}=\infty$ and
\begin{equation}\label{bound-inter}
\sup_{t\ge\tau}f(t)\le2f(\tau),\qquad\text{that is,}\qquad \|\nabla u(t)\|_{L^2}\le\sqrt2\,\e^{-(A(t)-A(\tau))}\|\nabla u(\tau)\|_{L^2},\quad t\ge\tau .
\end{equation}
\end{proposition}

\begin{proof}
We may assume that $u(\tau)\neq0$, since otherwise $u\equiv0$ on $[\tau,\infty)$ by uniqueness. As in the proof of Proposition \ref{prop:apriori-sub}, \eqref{basic}, \eqref{Nirenberg}, \eqref{eq:key} and the bound $\mathcal H(\tau)\le f(\tau)$ give, for $\tau\le t<T_{\max}$,
\begin{equation*}
F_\tau(t)\le f(\tau)+B_\tau F_\tau(t)^{\sigma/2},\qquad B_\tau=\frac{2K}{p+1}m^{p+1-\sigma}h(\tau).
\end{equation*}
We claim that $F_\tau(t)<2f(\tau)$ for all $t\in[\tau,T_{\max})$. The function $F_\tau$ is continuous and nondecreasing, and $F_\tau(\tau)=f(\tau)<2f(\tau)$. Note that $f(\tau)>0$, because a non-zero function of $\Ha$ cannot have a vanishing gradient. If the claim failed, there would be a first time $t_1>\tau$ with $F_\tau(t_1)=2f(\tau)$, and then
\begin{equation*}
2f(\tau)\le f(\tau)+B_\tau\big(2f(\tau)\big)^{\frac{\sigma-2}2}\,2f(\tau).
\end{equation*}
Now, we use $m=\e^{A(\tau)}\|u(\tau)\|_{L^2}$, $h(\tau)=\e^{-(p-1)A(\tau)}$ and $f(\tau)=\e^{2A(\tau)}\|\nabla u(\tau)\|_{L^2}^2$. Since $(p+1-\sigma)-(p-1)+(\sigma-2)=0$, we find
\begin{equation*}
B_\tau f(\tau)^{\frac{\sigma-2}2}=\frac{2K}{p+1}\|u(\tau)\|_{L^2}^{p+1-\sigma}\|\nabla u(\tau)\|_{L^2}^{\sigma-2}=\frac{2K}{p+1}\mathcal Q(u(\tau)).
\end{equation*}
It follows that $B_\tau(2f(\tau))^{\frac{\sigma-2}2}\,2f(\tau)=\frac{2K}{p+1}2^{\frac{\sigma-2}{2}}\mathcal Q(u(\tau))\,2f(\tau)=\frac{K\,2^{(\sigma+2)/2}}{p+1}\mathcal Q(u(\tau))\,f(\tau)<f(\tau)$, by the definition of $\varepsilon_0$. This contradicts the previous inequality and proves the claim. As before, the blow-up alternative gives $T_{\max}=\infty$, and \eqref{bound-inter} follows by letting $t\to\infty$.
\end{proof}

\begin{remark}\rm 
For $\tau=0$, Proposition \ref{prop:apriori-inter} is the global existence part of \cite[Theorem 2.3]{HM-MJM2025}. There, it was obtained for general nonnegative dampings by a bootstrap argument (\cite[Lemma 3.2]{HM-MJM2025}, after Strauss \cite{Strauss1968}). The point of allowing an arbitrary initial time $\tau$ is that the quantity $\mathcal Q(u(\tau))$ may become small for large $\tau$, because the mass $\|u(\tau)\|_{L^2}=m\e^{-A(\tau)}$ decays. This is precisely what drives Theorem \ref{thm:dichotomy}.
\end{remark}

\section{Proof of the scattering criterion}\label{S5}

\begin{proof}[Proof of Theorem \ref{thm:criterion}]
We write $(q,r,k)=(q_s,r_s,k_s)$ and $F=h|v|^{p-1}v$. Since $|\nabla(|v|^{p-1}v)|\le p|v|^{p-1}|\nabla v|$, H\"older's inequality in space with the exponents given by \eqref{Holder-rs} yields, for every $t\ge0$,
\begin{equation}\label{nonlin-est}
\|F(t)\|_{W^{1,r'}}\le p\,h(t)\|v(t)\|_{L^s}^{p-1}\|v(t)\|_{W^{1,r}}=p\,\Phi(t)\|v(t)\|_{W^{1,r}} .
\end{equation}
In case (a), H\"older's inequality in time with $\frac1{q'}=\frac1k+\frac1q$ gives, for every interval $J\subset[0,\infty)$,
\begin{equation}\label{nonlin-est2}
\|F\|_{L^{q'}(J;W^{1,r'})}\le p\|\Phi\|_{L^k(J)}\|v\|_{L^q(J;W^{1,r})}.
\end{equation}
In case (b), we have $k=\infty$ and $q=2$, and the same inequality holds with $\|\Phi\|_{L^\infty(J)}$.

\emph{Step 1: global Strichartz bounds in case (a).} Since $k<\infty$, we can split $[0,\infty)$ into consecutive intervals $I_j=[t_{j-1},t_j)$, $1\le j\le n'$, with $t_0=0$ and $t_{n'}=\infty$, such that $\|\Phi\|_{L^k(I_j)}^k=(2pC_{\rm St})^{-k}$ for $j<n'$ and $\|\Phi\|_{L^k(I_{n'})}^k\le(2pC_{\rm St})^{-k}$. Indeed, $t\mapsto\int_{t_{j-1}}^t\Phi^k$ is continuous and nondecreasing, so that each $t_j$ can be chosen as the first time at which this integral reaches $(2pC_{\rm St})^{-k}$, if such a time exists. Each of the first $n'-1$ intervals carries the amount $(2pC_{\rm St})^{-k}$ of $\|\Phi\|_{L^k(0,\infty)}^k$. We infer that $(n'-1)(2pC_{\rm St})^{-k}\le\|\Phi\|_{L^k(0,\infty)}^k$, so that the procedure stops after finitely many steps and $n'\le n$. For $1\le j\le n'$, we set $Y_j=\|v\|_{L^q(I_j;W^{1,r})}+\|v\|_{L^\infty(I_j;H^1)}$, and we fix $T\in I_j$. We apply \eqref{St1}, \eqref{St2} and the unitarity of $\e^{it\Delta}$ on $\Ha$ to the Duhamel formula \eqref{Integ-Eq} with $\tau=t_{j-1}$, and we use \eqref{nonlin-est2} with $J=[t_{j-1},T]$. We obtain
\begin{equation*}
\|v\|_{L^q(t_{j-1},T;W^{1,r})}+\|v\|_{L^\infty(t_{j-1},T;H^1)}\le 2C_{\rm St}\|v(t_{j-1})\|_{H^1}+\frac12\|v\|_{L^q(t_{j-1},T;W^{1,r})},
\end{equation*}
since $pC_{\rm St}\|\Phi\|_{L^k(I_j)}\le\frac12$. The left-hand side is finite by Proposition \ref{prop:lwp}. We deduce that it is bounded by $4C_{\rm St}\|v(t_{j-1})\|_{H^1}$, and letting $T\uparrow t_j$, we get $Y_j\le4C_{\rm St}\|v(t_{j-1})\|_{H^1}$. Since $\|v(t_j)\|_{H^1}\le Y_j$ for $j<n'$ by continuity of $v$, an induction gives $Y_j\le(4C_{\rm St})^j\|u_0\|_{H^1}$. Summing over $j$, we obtain \eqref{eq:Strbound}, since $(4C_{\rm St})^j\le(4C_{\rm St})^n$ and $n'\le n$. In particular, $v\in L^q(0,\infty;W^{1,r})$.

\emph{Step 2: global Strichartz bounds in case (b).} We choose $\tau_0\ge0$ such that $pC_{\rm St}\sup_{t\ge\tau_0}\Phi(t)\le\frac12$, which is possible since $\limsup\Phi<\delta_*$. Arguing as in Step 1 on the interval $[\tau_0,T]$, we obtain $\|v\|_{L^q(\tau_0,T;W^{1,r})}\le4C_{\rm St}\|v(\tau_0)\|_{H^1}$ for every $T>\tau_0$. Letting $T\to\infty$, we get $v\in L^q(\tau_0,\infty;W^{1,r})$. Since $v\in L^q(0,\tau_0;W^{1,r})$ by Proposition \ref{prop:lwp}, we conclude that $v\in L^q(0,\infty;W^{1,r})$.

\emph{Step 3: conclusion.} In both cases, \eqref{nonlin-est2} shows that $F\in L^{q'}(0,\infty;W^{1,r'})$. Applying \eqref{St1}--\eqref{St2} once more, with an arbitrary admissible pair on the left-hand side, we infer that $v\in L^{\tilde q}(0,\infty;W^{1,\tilde r})$ for every admissible pair $(\tilde q,\tilde r)$. We now prove scattering. By \eqref{Integ-Eq} with $\tau=0$, we have $\e^{-it\Delta}v(t)=u_0+i\int_0^t\e^{-is\Delta}F(s)\dd s$. Hence, for $0\le t_1<t_2$, \eqref{St3} and \eqref{nonlin-est2} give
\begin{equation*}
\big\|\e^{-it_2\Delta}v(t_2)-\e^{-it_1\Delta}v(t_1)\big\|_{H^1}=\Big\|\int_{t_1}^{t_2}\e^{-is\Delta}F(s)\dd s\Big\|_{H^1}\le pC_{\rm St}\|\Phi\|_{L^k(t_1,t_2)}\|v\|_{L^q(t_1,t_2;W^{1,r})}.
\end{equation*}
In both cases, the right-hand side tends to zero as $t_1,t_2\to\infty$. Indeed, $\|\Phi\|_{L^k(t_1,t_2)}$ is bounded, by $\|\Phi\|_{L^k(0,\infty)}$ in case (a) and by $\sup\Phi\le L_s^{p-1}$ in case (b), while $v\in L^q(0,\infty;W^{1,r})$ with $q<\infty$. Therefore, $\e^{-it\Delta}v(t)$ is a Cauchy family in $\Ha$ as $t\to\infty$, and it converges to some $v_+\in\Ha$. Letting $t_2\to\infty$ in the last inequality with $t_1=t$, we obtain \eqref{rate}.
\end{proof}

\begin{remark}\rm 
The proof of Theorem \ref{thm:criterion} is the time-dependent analogue of \cite[Lemmas 3.4 and 3.5]{HM-MJM2025}, where the role of $h$ was played by the exponential weight $\e^{-(p-1)\Lsub t}$. Here, the decay is placed on the coupling itself, which makes it possible to deal with dampings for which $\Lsub=0$. Compared with the Lorentz-space approach of \cite{BamriTayachi2023}, the criterion requires an a priori bound on $v$ in $L^\infty(0,\infty;L^s)$, but no smallness of the data. Note also that, in case (a), the bound \eqref{eq:Strbound} depends only on $\|u_0\|_{H^1}$ and on $\|\Phi\|_{L^{k_s}(0,\infty)}$, and not on the local theory.
\end{remark}

\section{Proof of the non-scattering criterion}\label{S6}

\begin{proof}[Proof of Theorem \ref{thm:nonscat}]
Recall that $\sigma=\frac{N(p-1)}2$, and let $\varrho=\frac2{2-p}\in(2,\infty]$, so that $\varrho'=\frac2p\in[1,2)$ and $N\big(\frac12-\frac1\varrho\big)=\sigma$. By the conservation of mass for \eqref{V}, we have $\|v(t)\|_{L^2}=\|u_0\|_{L^2}=m$ for all $t$. We infer that $\|v_+\|_{L^2}=\lim_{t\to\infty}\|\e^{-it\Delta}v(t)\|_{L^2}=m$. We argue by contradiction and assume that $m>0$.

Set $g=|\hat v_+|^{p-1}\hat v_+$. Since $\hat v_+\in\Lt$, we have $g\in L^{2/p}(\R^N)\subset L^1_{\rm loc}(\R^N)$ and $g\neq0$. Hence there exists $\psi\in\cS(\R^N)$ such that
\begin{equation}\label{choice-psi}
\int_{\R^N}g(\xi)\,\overline{\hat\psi(\xi)}\dd\xi=1.
\end{equation}
Indeed, otherwise $g$ would vanish as a distribution, since $\Fou$ is an isomorphism of $\cS(\R^N)$. Consider
\begin{equation*}
\phi(t)=\big(\e^{-it\Delta}v(t),\psi\big)=\big(v(t),\e^{it\Delta}\psi\big),\qquad t\ge0 .
\end{equation*}
By assumption, $\phi(t)\to(v_+,\psi)$ as $t\to\infty$. In particular, $\phi$ is bounded. We now compute the derivative of $\phi$. Since $p<1+\frac4{(N-2)_+}$, the nonlinearity maps $\Ha$ continuously into $H^{-1}(\R^N)$, and the solution satisfies $v\in C^1([0,\infty);H^{-1}(\R^N))$ with $\partial_tv=i\Delta v+ih|v|^{p-1}v$ in $H^{-1}$; see \cite[Section 4.4]{Cazenave2003}. Moreover, $w(t)=\e^{it\Delta}\psi$ belongs to $C^1([0,\infty);\Ha)$, with $\partial_tw=i\Delta w$. Hence $\phi$ is of class $C^1$, and
\begin{equation*}
\phi'(t)=\big(i\Delta v+ih|v|^{p-1}v,w\big)+\big(v,i\Delta w\big)=ih(t)\big(|v|^{p-1}v,w\big),
\end{equation*}
because $(i\Delta v,w)+(v,i\Delta w)=i(\Delta v,w)-i(v,\Delta w)=0$. Integrating, we obtain
\begin{equation}\label{phi-deriv}
\phi(t_2)-\phi(t_1)=i\int_{t_1}^{t_2}h(t)\,\mathcal N(t)\dd t,\qquad \mathcal N(t)=\big(|v(t)|^{p-1}v(t),\e^{it\Delta}\psi\big).
\end{equation}
The function $\mathcal N$ is continuous on $[0,\infty)$. Indeed, $v\in C([0,\infty);\Lt)$, the map $w\mapsto|w|^{p-1}w$ is continuous from $\Lt$ into $L^{2/p}$, and $t\mapsto\e^{it\Delta}\psi$ is continuous with values in $L^{\varrho}$.

We now compute the asymptotic behaviour of $\mathcal N(t)$. Let $w(t)=\e^{it\Delta}\psi$ and $Z(t)=\Mo(t)\Do(t)\hat v_+$. By \eqref{MDFM}, we have $\e^{it\Delta}v_+-Z(t)=\Mo(t)\Do(t)\Fou\big((\Mo(t)-1)v_+\big)$, so that $\|\e^{it\Delta}v_+-Z(t)\|_{L^2}=\|(\Mo(t)-1)v_+\|_{L^2}\to0$ by dominated convergence. Consequently, $\rho(t)=v(t)-Z(t)$ satisfies
\begin{equation*}
\|\rho(t)\|_{L^2}\le\|\e^{-it\Delta}v(t)-v_+\|_{L^2}+\|(\Mo(t)-1)v_+\|_{L^2}\xrightarrow[t\to\infty]{}0 .
\end{equation*}
Similarly, $w(t)-\Mo(t)\Do(t)\hat\psi=\Mo(t)\Do(t)\Fou\big((\Mo(t)-1)\psi\big)$. We use \eqref{MDFM}, the Hausdorff--Young inequality (whose constant is at most $1$ for the unitary normalisation of $\Fou$ and $\varrho\ge2$) and the bound $|\e^{i|x|^2/(4t)}-1|\le\frac{|x|^2}{4t}$. We get
\begin{equation*}
\big\|w(t)-\Mo(t)\Do(t)\hat\psi\big\|_{L^\varrho}\le(2t)^{-\sigma}\big\|(\Mo(t)-1)\psi\big\|_{L^{\varrho'}}\le\frac{(2t)^{-\sigma}}{4t}\,\big\||x|^2\psi\big\|_{L^{\varrho'}} .
\end{equation*}
Moreover, \eqref{dispersive} gives $\|w(t)\|_{L^\varrho}\le(4\pi t)^{-\sigma}\|\psi\|_{L^{\varrho'}}$. We decompose
\begin{equation*}
\mathcal N(t)=\big(|v|^{p-1}v-|Z|^{p-1}Z,\,w\big)+\big(|Z|^{p-1}Z,\,w-\Mo\Do\hat\psi\big)+\big(|Z|^{p-1}Z,\,\Mo\Do\hat\psi\big)=:\mathcal N_1+\mathcal N_2+\mathcal N_3 .
\end{equation*}
We use the elementary inequality $\big||z_1|^{p-1}z_1-|z_2|^{p-1}z_2\big|\le p\big(|z_1|^{p-1}+|z_2|^{p-1}\big)|z_1-z_2|$ and H\"older's inequality with $\frac{p-1}2+\frac12+\frac1\varrho=1$, together with $\|v(t)\|_{L^2}=\|Z(t)\|_{L^2}=m$. We obtain
\begin{equation*}
|\mathcal N_1(t)|\le 2p\,m^{p-1}\|\rho(t)\|_{L^2}\|w(t)\|_{L^\varrho}\le C\,t^{-\sigma}\|\rho(t)\|_{L^2}=o(t^{-\sigma}).
\end{equation*}
Since $\||Z|^{p-1}Z\|_{L^{2/p}}=m^p$, we also get
\begin{equation*}
|\mathcal N_2(t)|\le m^p\,\big\|w(t)-\Mo(t)\Do(t)\hat\psi\big\|_{L^\varrho}\le C\,t^{-\sigma-1}.
\end{equation*}
Finally, the phases $\Mo(t)$ cancel in $\mathcal N_3$, as does the constant phase of $(2i)^{-N/2}$ in $\Do(t)$. The change of variables $x=2t\xi$ then gives
\begin{equation*}
\mathcal N_3(t)=\int_{\R^N}|\Do(t)\hat v_+|^{p-1}\Do(t)\hat v_+\,\overline{\Do(t)\hat\psi}\dd x=(2t)^{-\sigma}\int_{\R^N}g(\xi)\overline{\hat\psi(\xi)}\dd\xi=(2t)^{-\sigma},
\end{equation*}
by \eqref{choice-psi}. Hence $\mathcal N(t)=(2t)^{-\sigma}\big(1+\epsilon(t)\big)$ with $\epsilon(t)\to0$, and there exists $T_1\ge1$ such that ${\rm Re}\,\mathcal N(t)\ge\frac12(2t)^{-\sigma}$ for $t\ge T_1$. Taking imaginary parts in \eqref{phi-deriv} with $t_1=T_1$, and using $h\ge0$, we find
\begin{equation*}
{\rm Im}\,\phi(t_2)-{\rm Im}\,\phi(T_1)=\int_{T_1}^{t_2}h(t)\,{\rm Re}\,\mathcal N(t)\dd t\ge2^{-\sigma-1}\int_{T_1}^{t_2}h(t)\,t^{-\sigma}\dd t\xrightarrow[t_2\to\infty]{}\infty ,
\end{equation*}
by \eqref{nonint}, since $h(t)t^{-\sigma}$ is integrable on $[1,T_1]$. This contradicts the boundedness of $\phi$. Therefore $m=0$, that is, $v\equiv0$.
\end{proof}

\begin{remark}\rm 
The restriction $p\le2$ is used twice. First, it ensures $\varrho=\frac2{2-p}\ge2$, so that the dispersive estimate and the Hausdorff--Young inequality apply to $w$. Second, it is needed in the H\"older inequality for $\mathcal N_1$, which only uses $L^2$ information on $v-Z$. It is the same restriction ($\alpha=p-1\le1$) as in \cite[Theorem 1.7 and Lemma 6.2]{BamriTayachi2023}. For $N\ge2$, it is harmless in our applications, since $1+\frac2N\le2$. For $N=1$, it excludes the range $2<p\le3$, which is covered by Barab's argument \cite{Barab1984} when $h\equiv1$.
\end{remark}

\section{Proofs of the main results}\label{S7}

Throughout this section, $\da(t)=a(1+t)^{-\theta}$, $u$ is the maximal solution of \eqref{Main-New}, $v=\e^{A}u$ solves \eqref{V} with $h=\e^{-(p-1)A}$, and $m=\|u_0\|_{L^2}$. We shall repeatedly use the following consequence of Theorem \ref{thm:criterion}.

\begin{lemma}\label{lem:bounded-scatter}
Assume $(\mathrm D_{N,p})$, and let $u$ be a global solution of \eqref{Main-New} such that $V:=\sup_{t\ge0}\e^{A(t)}\|u(t)\|_{H^1}<\infty$. Then $u$ scatters in $\Ha$. Moreover, if $s$ is the exponent given by Lemma \ref{lem:Dimplies}, then
\begin{equation}\label{rate2}
\big\|\e^{A(t)}\e^{-it\Delta}u(t)-u_+\big\|_{H^1}\le C\,\|h\|_{L^{k_s}(t,\infty)},\qquad t\ge0,
\end{equation}
where $C$ depends only on $N$, $p$, $s$, $V$ and $\|h\|_{L^{k_s}(0,\infty)}$.
\end{lemma}

\begin{proof}
By assumption, $\sup_{t\ge0}\|v(t)\|_{H^1}=V<\infty$. By the Sobolev embedding $\Ha\hookrightarrow L^s$, valid for $s\in J_N$, we get $L_s=\sup_t\|v(t)\|_{L^s}\le C_sV$. Hence $\Phi(t)=h(t)\|v(t)\|_{L^s}^{p-1}\le L_s^{p-1}h(t)$ belongs to $L^{k_s}(0,\infty)$ by Lemma \ref{lem:Dimplies}, and case (a) of Theorem \ref{thm:criterion} applies. Estimate \eqref{rate2} follows from \eqref{rate} and \eqref{eq:Strbound}, since $\e^{A(t)}\e^{-it\Delta}u(t)=\e^{-it\Delta}v(t)$, $\|\Phi\|_{L^{k_s}(t,\infty)}\le L_s^{p-1}\|h\|_{L^{k_s}(t,\infty)}$ and $\|u_0\|_{H^1}\le V$, the integer $n$ in \eqref{eq:Strbound} being a nondecreasing function of $\|\Phi\|_{L^{k_s}(0,\infty)}\le L_s^{p-1}\|h\|_{L^{k_s}(0,\infty)}$.
\end{proof}

\begin{proof}[Proof of Theorem \ref{thm:sub}]
Part (i) is Proposition \ref{prop:apriori-sub} together with \eqref{M-Id}, since $\e^{A(t)}\|\nabla u(t)\|_{L^2}=\|\nabla v(t)\|_{L^2}$. Part (ii) follows from part (i) and Lemma \ref{lem:bounded-scatter}, the explicit form of $(\mathrm D_{N,p})$ being given in Remark \ref{rem:Dexplicit}. For the rate \eqref{rate-sub}, let $\theta=1$ and $p>1+\frac4{N+4a}$, with $N\ge2$, or with $N=1$ and $p\le3$. We choose $s=2$. This exponent is $p$-admissible with $q_2>2$ and $\frac1{k_2}=1-\frac{N(p-1)}4$ (see the proof of Lemma \ref{lem:ell}), and the condition $p>1+\frac4{N+4a}$ is exactly $k_2a(p-1)>1$. Then $\Phi=m^{p-1}h$, and \eqref{rate}, \eqref{tail-h} give
\begin{equation*}
\big\|\e^{A(t)}\e^{-it\Delta}u(t)-u_+\big\|_{H^1}\le pC_{\rm St}m^{p-1}\big(k_2a(p-1)-1\big)^{-\frac1{k_2}}\|v\|_{L^{q_2}(0,\infty;W^{1,r_2})}\,(1+t)^{-a(p-1)+1-\frac{N(p-1)}4},
\end{equation*}
which is \eqref{rate-sub}. By \eqref{eq:Strbound}, the Strichartz norm on the right-hand side is bounded by $2n(4C_{\rm St})^n\|u_0\|_{H^1}$, where $n$ only depends on $N,p,a$ and $m$, through $\|\Phi\|_{L^{k_2}(0,\infty)}=m^{p-1}\big(k_2a(p-1)-1\big)^{-1/k_2}$. Hence the constant only depends on $N,p,a$ and $\|u_0\|_{H^1}$.

For part (iii), assume that $u_0\ne0$, $p\le2$, and that $u$ scatters in $L^2$, that is, $\|\e^{-it\Delta}v(t)-u_+\|_{L^2}\to0$ for some $u_+\in\Lt$. The solution $v$ is global by part (i), and Theorem \ref{thm:nonscat} yields $v\equiv0$, which is a contradiction, provided that \eqref{nonint} holds. If $\theta=1$, then $2^{-a(p-1)}t^{-(p-1)(a+\frac N2)}\le h(t)t^{-\frac{N(p-1)}2}\le t^{-(p-1)(a+\frac N2)}$ for $t\ge1$. We deduce that \eqref{nonint} holds if and only if $(p-1)(a+\frac N2)\le1$, that is, $p\le1+\frac2{N+2a}$. If $\theta>1$, then $h\ge h_\infty>0$ by Lemma \ref{lem:damping}(iii), and \eqref{nonint} holds as soon as $\frac{N(p-1)}2\le1$, that is, $p\le1+\frac2N$. In both cases $p<1+\frac4N$, so that part (i) indeed applies.
\end{proof}

\begin{proof}[Proof of Theorem \ref{thm:crit}]
Part (i) is Proposition \ref{prop:apriori-crit} with $\tau=0$, for which $\eta_0=(m/\|Q\|_{L^2})^{4/N}$ and $\mathcal H(0)=E(u_0)$. For part (ii), assume that $T_{\max}>T_Q$, and let $\tau\in(T_Q,T_{\max})$. By the definition of $T_Q$, there exists $t'\in[T_Q,\tau)$ such that $\e^{-A(t')}m<\|Q\|_{L^2}$. Since $t\mapsto\e^{-A(t)}m$ is nonincreasing, we infer that $\|u(\tau)\|_{L^2}=\e^{-A(\tau)}m<\|Q\|_{L^2}$. Proposition \ref{prop:apriori-crit} then gives $T_{\max}=\infty$ and $\sup_{t\ge\tau}\|\nabla v(t)\|_{L^2}<\infty$, while $\|\nabla v(t)\|_{L^2}$ is bounded on $[0,\tau]$ by continuity. Together with $\|v(t)\|_{L^2}=m$, we obtain $\sup_{t\ge0}\e^{A(t)}\|u(t)\|_{H^1}<\infty$. If $\theta\le1$ and $a>0$, then $T_Q<\infty$ by \eqref{eq:TQ-explicit}, so that every global solution satisfies $T_{\max}=\infty>T_Q$. Part (iii) follows from part (ii) and Lemma \ref{lem:bounded-scatter}, the explicit form of $(\mathrm D_{N,p})$ for $p=1+\frac4N$ being given in Remark \ref{rem:Dexplicit}.
\end{proof}

\begin{proof}[Proof of Theorem \ref{thm:inter-small}]
Global existence and \eqref{Decay-u} follow from Proposition \ref{prop:apriori-inter} with $\tau=0$ and from \eqref{M-Id}. In case (a), scattering follows from Lemma \ref{lem:bounded-scatter}, since \eqref{Decay-u} gives $\sup_{t\ge0}\e^{A(t)}\|u(t)\|_{H^1}\le\sqrt2\,\|u_0\|_{H^1}$.

In case (b), let $N\ge3$ and $s_*=\frac{N(p-1)}2\in(2,\frac{2N}{N-2})$. Then $s_*$ is $p$-admissible, $(q_{s_*},r_{s_*})=(2,\frac{2N}{N-2})$ is the endpoint pair, and $k_{s_*}=\infty$. We use \eqref{eq:GNcrit}, \eqref{Decay-u} in the form $\|\nabla v(t)\|_{L^2}\le\sqrt2\|\nabla u_0\|_{L^2}$, the bound $h\le1$, and the identities $(1-s_c)(p-1)=p+1-\sigma$ and $s_c(p-1)=\sigma-2$. We obtain, for all $t\ge0$,
\begin{equation*}
\Phi(t)=h(t)\|v(t)\|_{L^{s_*}}^{p-1}\le K_*^{p-1}m^{(1-s_c)(p-1)}\big(\sqrt2\|\nabla u_0\|_{L^2}\big)^{s_c(p-1)}=K_*^{p-1}2^{\frac{\sigma-2}2}\mathcal Q(u_0).
\end{equation*}
We set $\varepsilon_1=\min\big\{\varepsilon_0,\ \frac12\delta_*K_*^{1-p}2^{-\frac{\sigma-2}2}\big\}$, where $\delta_*$ is the constant of Theorem \ref{thm:criterion}(b). If $\mathcal Q(u_0)<\varepsilon_1$, then $\sup_t\Phi(t)<\delta_*$, and $L_{s_*}<\infty$ by the Sobolev embedding. Therefore, case (b) of Theorem \ref{thm:criterion} gives scattering. No condition on $a\ge0$ or $\theta\ge0$ has been used.
\end{proof}

\begin{proof}[Proof of Theorem \ref{thm:dichotomy}]
Assume first that $\mathcal Q(u(\tau))<\varepsilon_0$ for some $\tau\ge0$. By Proposition \ref{prop:apriori-inter}, we have $\sup_{t\ge\tau}\|\nabla v(t)\|_{L^2}^2\le2\|\nabla v(\tau)\|_{L^2}^2$. Since $\|\nabla v\|_{L^2}$ is bounded on $[0,\tau]$ by continuity and $\|v(t)\|_{L^2}=m$, we obtain $\sup_{t\ge0}\e^{A(t)}\|u(t)\|_{H^1}<\infty$. Lemma \ref{lem:bounded-scatter} then shows that $u$ scatters, so that alternative (a) holds.

Otherwise, $\mathcal Q(u(t))\ge\varepsilon_0$ for every $t\ge0$. Since $\|u(t)\|_{L^2}=m\e^{-A(t)}$ and $p+1-\sigma>0$, this reads
\begin{equation*}
\|\nabla u(t)\|_{L^2}^{\sigma-2}\ge\varepsilon_0\,\|u(t)\|_{L^2}^{-(p+1-\sigma)}=\varepsilon_0\,m^{-(p+1-\sigma)}\e^{(p+1-\sigma)A(t)},
\end{equation*}
which is \eqref{growup} after taking the power $\frac1{\sigma-2}$. Since $A(t)\to\infty$ by Lemma \ref{lem:Dimplies}, we get $\|\nabla u(t)\|_{L^2}\to\infty$. By Lemma \ref{lem:liminf}, there exists a sequence $t_n\to\infty$ such that $I(u(t_n))\le\frac1n$, that is, $\|u(t_n)\|_{L^{p+1}}^{p+1}\ge\|\nabla u(t_n)\|_{L^2}^2-\frac1n\to\infty$. Hence $\limsup_{t\to\infty}\|u(t)\|_{L^{p+1}}=\infty$.

The two alternatives are mutually exclusive. Indeed, if $u$ scatters, then $\|\nabla u(t)\|_{L^2}\le\e^{-A(t)}\sup_s\|v(s)\|_{H^1}\to0$ by Remark \ref{rem:def}(i), which is incompatible with \eqref{growup}.
\end{proof}

\begin{proof}[Proof of Corollary \ref{cor:Lp+1}]
Each of the three conditions in \eqref{hyp-sup} excludes alternative (b) of Theorem \ref{thm:dichotomy}. The first one does, because (b) implies $\limsup\|u(t)\|_{L^{p+1}}=\infty$. The second one does, because of \eqref{growup}. The third one implies the first one, since $\Ha\hookrightarrow L^{p+1}$. Hence alternative (a) holds.
\end{proof}

\begin{proof}[Proof of Corollary \ref{cor:general}]

For (i), we choose $a'$ with $\frac{\ell_N(p)}{p-1}<a'<\liminf_{t\to\infty}t\,\da(t)$, and $T\ge1$ such that $t\,\da(t)\ge a'$ for $t\ge T$. Then $A(t)\ge A(T)+a'\log(t/T)$ for $t\ge T$, and we obtain $h(t)\le C\,t^{-a'(p-1)}$ for $t\ge T$. By Lemma \ref{lem:ell}, there exists a $p$-admissible $s$ with $q_s>2$ and $\frac1{k_s}<a'(p-1)$. Since $h\le1$, we deduce that $h\in L^{k_s}(0,\infty)$. The conclusion follows from Remark \ref{rem:general}. Indeed, the a priori bounds of Section \ref{S4} hold for every nonnegative damping, and the proofs given above apply without change; in particular, $A(t)\to\infty$, so that $T_Q<\infty$ in the mass-critical case and Lemma \ref{lem:liminf} applies in the intercritical case. For (ii), we have $A(t)\le A(T)+b\log(t/T)$ for $t\ge T$, so that $h(t)\ge c\,t^{-b(p-1)}$ with $c>0$. Since $(p-1)(b+\frac N2)\le1$, the integral \eqref{nonint} diverges. Moreover, the condition $(p-1)\frac N2\le1$ gives $p\le1+\frac2N<1+\frac4N$. Hence every solution is global by Proposition \ref{prop:apriori-sub}, and Theorem \ref{thm:nonscat} applies.
\end{proof}

\section{The scale-invariant case: structure of the gap}\label{S-gap}

In this section, we take a closer look at the range left open by Theorem \ref{thm:sub} in the scale-invariant case. Throughout, $\theta=1$, $a>0$ and $1<p<1+\frac4N$, so that $h(t)=(1+t)^{-\alpha}$ with
\begin{equation}\label{eq:alphabeta}
\alpha=a(p-1),\qquad \sigma=\frac{N(p-1)}2<2,\qquad \beta=\frac{\alpha}{2-\sigma}=\frac{2a(p-1)}{4-N(p-1)}>0 .
\end{equation}
Elementary computations show that
\begin{equation}\label{eq:equivalences}
\begin{gathered}
p\le p_{\rm ns}(a)\iff\alpha+\sigma\le1,\qquad p<p_{\rm s}(a)\iff\alpha+\frac\sigma2<1\iff\beta<\frac12,\\
p<1+\frac4{N+2a}\iff\beta<1 .
\end{gathered}
\end{equation}
We first explain where the two thresholds come from. We then reduce the problem exactly to a nonlinear Schr\"odinger equation with a repulsive harmonic potential, and we derive some rigorous consequences. Finally, we formulate a precise conjecture on the gap and present numerical evidence for it. We insist on the separation between the proved statements (Lemma \ref{lem:scaling}, Propositions \ref{prop:reduction}--\ref{prop:approx}) and the conjectural ones (Conjecture \ref{conj:adiabatic}).

\subsection{Scaling and the meaning of the two thresholds}
The coupling $h(t)=(1+t)^{-\alpha}$ is asymptotic to the homogeneous coupling $t^{-\alpha}$, for which \eqref{V} enjoys a scaling invariance.

\begin{lemma}\label{lem:scaling}
Let $\alpha\ge0$ and let $v$ be a solution of $i\partial_tv+\Delta v+t^{-\alpha}|v|^{p-1}v=0$ on $(0,\infty)\times\R^N$. Then, for every $\mu>0$, the function $v_\mu(t,x)=\mu^{\frac{2(1-\alpha)}{p-1}}v(\mu^2t,\mu x)$ is a solution of the same equation, and
\begin{equation*}
\|v_\mu(t)\|_{L^2}=\mu^{\frac{2(1-\alpha)}{p-1}-\frac N2}\,\|v(\mu^2t)\|_{L^2}.
\end{equation*}
In particular, the scaling is $L^2$-invariant if and only if $\alpha+\frac\sigma2=1$, that is, when $\alpha=a(p-1)$, if and only if $p=p_{\rm s}(a)$.
\end{lemma}

\begin{proof}
Set $\kappa'=\frac{2(1-\alpha)}{p-1}$. Then $i\partial_tv_\mu+\Delta v_\mu=\mu^{\kappa'+2}\big(i\partial_tv+\Delta v\big)(\mu^2t,\mu x)$, while $t^{-\alpha}|v_\mu|^{p-1}v_\mu=\mu^{\kappa'p+2\alpha}\big((\mu^2t)^{-\alpha}|v|^{p-1}v\big)(\mu^2t,\mu x)$. Since $\kappa'p+2\alpha=\kappa'+2$, the first claim follows. The formula for the $L^2$ norm is a change of variables, and the exponent $\kappa'-\frac N2$ vanishes if and only if $2(1-\alpha)=\sigma$.
\end{proof}

Lemma \ref{lem:scaling} identifies $p_{\rm s}(a)$ as the \emph{mass-critical exponent} of the scale-invariant problem, exactly as $1+\frac4N$ is the mass-critical exponent of the undamped equation. Similarly, along a free evolution $\e^{it\Delta}\varphi$, the effective potential $h(t)|\e^{it\Delta}\varphi|^{p-1}$ has size $t^{-(\alpha+\sigma)}$, which is integrable at infinity if and only if $p>p_{\rm ns}(a)$. Thus $p_{\rm ns}(a)$ is the \emph{short-range threshold}, exactly as $1+\frac2N$ for the undamped equation. In other words, the scale-invariant damping shifts both classical thresholds by replacing $\sigma$ with $\sigma+\alpha$ in the long-range condition and with $\frac\sigma2+\alpha$ in the mass-critical condition. Table \ref{tab:analogy} summarises this analogy. In the undamped focusing case, the range between the two thresholds is precisely where solitary waves prevent the scattering of all solutions, while small data scatter. The results below indicate that the same structure persists in the scale-invariant case, with solitary waves replaced by \emph{adiabatically spreading} solitary waves.

\begin{table}[ht]
\centering
\renewcommand{\arraystretch}{1.45}
{\small
\begin{tabular}{|p{4.6cm}|p{4.1cm}|p{4.9cm}|}
\hline
 & undamped focusing NLS ($a=0$) & scale-invariant damping ($\theta=1$, $a>0$)\\
\hline
coupling $h(t)$ & $1$ & $(1+t)^{-a(p-1)}$\\
\hline
short-range threshold (no free solution at or below) & $1+\frac2N$ ($\sigma=1$) & $p_{\rm ns}(a)=1+\frac2{N+2a}$ ($\alpha+\sigma=1$), Theorem \ref{thm:sub}(iii)\\
\hline
small data scattering & above the Strauss exponent & above $1+\alpha_1(2a)$ \cite{BamriTayachi2023}, see \eqref{eq:quadratic}\\
\hline
$L^2$-critical exponent of the scaling & $1+\frac4N$ ($\sigma=2$) & $p_{\rm s}(a)=1+\frac4{N+4a}$ ($\alpha+\frac\sigma2=1$), Lemma \ref{lem:scaling}\\
\hline
all solutions scatter & never (focusing) & $p>p_{\rm s}(a)$, Theorem \ref{thm:sub}(ii)\\
\hline
obstruction between the thresholds & solitary waves $\e^{i\omega t}R_\omega$ & adiabatic solitary waves spreading like $(1+t)^\beta$ (Conjecture \ref{conj:adiabatic})\\
\hline
pseudo-conformal energy nonincreasing & $p\le1+\frac4N$ & $p\le1+\frac4{N+2a}$, Proposition \ref{prop:pseudoconf}\\
\hline
\end{tabular}}
\caption{The scale-invariant damping shifts the classical exponents of the focusing NLS equation. The comparison suggests that the gap $(p_{\rm ns}(a),p_{\rm s}(a)]$ plays the role of the range $(1+\frac2N,1+\frac4N)$ in the undamped case.}
\label{tab:analogy}
\end{table}

\subsection{An exact reduction}
We now show that \eqref{Main-New} with $\theta=1$ is exactly equivalent to a nonlinear Schr\"odinger equation with a time-dependent repulsive harmonic potential. Let
\begin{equation}\label{eq:lambdatau}
\lambda(t)=(1+t)^{-\beta},\qquad \tau(t)=\int_0^t\lambda(s)^2\dd s=\begin{cases}\dfrac{(1+t)^{1-2\beta}-1}{1-2\beta}, &\beta\neq\frac12,\\[8pt] \log(1+t), &\beta=\frac12,\end{cases}
\end{equation}
and let $\tau_\infty=\lim_{t\to\infty}\tau(t)$. Thus $\tau_\infty=\infty$ if $\beta\le\frac12$, and $\tau_\infty=\frac1{2\beta-1}$ if $\beta>\frac12$. Note that $1+(1-2\beta)\tau(t)=(1+t)^{1-2\beta}$ for every $\beta$.

\begin{proposition}\label{prop:reduction}
Let $\theta=1$, $a>0$ and $1<p<1+\frac4N$, and let $\alpha,\sigma,\beta$ be given by \eqref{eq:alphabeta}. Let $u\in C([0,T);\Ha)$, $v=\e^{A}u$, and define $\widetilde Z$ on $[0,\tau(T))\times\R^N$ by
\begin{equation}\label{eq:Ztilde}
v(t,x)=\lambda(t)^{\frac N2}\exp\Big(\frac{i\beta|x|^2}{4(1+t)}\Big)\,\widetilde Z\big(\tau(t),\lambda(t)x\big).
\end{equation}
Then $u$ solves \eqref{Main-New} on $[0,T)$ if and only if $\widetilde Z$ solves, in the sense of distributions on $(0,\tau(T))\times\R^N$,
\begin{equation}\label{eq:repulsive}
i\partial_\tau\widetilde Z+\Delta\widetilde Z+|\widetilde Z|^{p-1}\widetilde Z+q(\tau)|y|^2\widetilde Z=0,\qquad q(\tau)=\frac{\beta(1-\beta)}{4\big(1+(1-2\beta)\tau\big)^2}.
\end{equation}
Moreover, $\|\widetilde Z(\tau(t))\|_{L^2}=\|v(t)\|_{L^2}=m$ for all $t$. More generally, if $v$ solves $i\partial_tv+\Delta v+h|v|^{p-1}v=\mathfrak r$ with a source term $\mathfrak r$, then \eqref{eq:repulsive} holds with right-hand side $F$, where $\mathfrak r(t,x)=\lambda^{\frac N2+2}\exp\big(\frac{i\beta|x|^2}{4(1+t)}\big)F(\tau(t),\lambda(t) x)$.
\end{proposition}

\begin{proof}
We proceed in two steps, and we write $\Lambda=\frac N2+y\cdot\nabla$ for the generator of $L^2$-dilations.

\emph{Step 1: rescaling.} Let $Z$ be defined by $v(t,x)=\lambda(t)^{N/2}Z(\tau(t),\lambda(t)x)$. A direct computation gives
\begin{equation*}
i\partial_tv+\Delta v+h|v|^{p-1}v=\lambda^{\frac N2}\Big[i\tau'\partial_\tau Z+i\frac{\lambda'}{\lambda}\Lambda Z+\lambda^2\Delta Z+h\lambda^{\sigma}|Z|^{p-1}Z\Big](\tau,\lambda x),
\end{equation*}
where we used $|v|^{p-1}=\lambda^{N(p-1)/2}|Z|^{p-1}=\lambda^\sigma|Z|^{p-1}$. Now $\tau'=\lambda^2$ by definition, and $h\lambda^{\sigma}=\lambda^2$ because $\lambda^{2-\sigma}=(1+t)^{-\beta(2-\sigma)}=(1+t)^{-\alpha}=h$. We obtain
\begin{equation*}
i\partial_tv+\Delta v+h|v|^{p-1}v=\lambda^{\frac N2+2}\Big[i\partial_\tau Z+\Delta Z+|Z|^{p-1}Z+i\varepsilon\Lambda Z\Big](\tau,\lambda x),
\end{equation*}
with
\begin{equation*}
\varepsilon=\frac{\lambda'}{\lambda^3}=-\beta(1+t)^{2\beta-1}=-\frac{\beta}{1+(1-2\beta)\tau}.
\end{equation*}

\emph{Step 2: removal of the dilation term.} We write $Z=\e^{-ib|y|^2}\widetilde Z$ with $b=\frac\varepsilon4$, which is a function of $\tau$. Using $\nabla(\e^{-ib|y|^2})=-2iby\,\e^{-ib|y|^2}$, we compute
\begin{equation*}
\begin{aligned}
i\partial_\tau Z&=\e^{-ib|y|^2}\big(i\partial_\tau\widetilde Z+\partial_\tau b\,|y|^2\widetilde Z\big),\\
\Delta Z&=\e^{-ib|y|^2}\big(\Delta\widetilde Z-4iby\cdot\nabla\widetilde Z-2ibN\widetilde Z-4b^2|y|^2\widetilde Z\big),\\
i\varepsilon\Lambda Z&=\e^{-ib|y|^2}\big(i\varepsilon\Lambda\widetilde Z+2\varepsilon b|y|^2\widetilde Z\big).
\end{aligned}
\end{equation*}
Since $\varepsilon=4b$, the terms involving $y\cdot\nabla\widetilde Z$ and $N\widetilde Z$ cancel. The coefficient of $|y|^2\widetilde Z$ is $\partial_\tau b-4b^2+2\varepsilon b=\frac14(\partial_\tau\varepsilon+\varepsilon^2)$. Finally,
\begin{equation*}
\partial_\tau\varepsilon+\varepsilon^2=\frac{\beta(1-2\beta)+\beta^2}{\big(1+(1-2\beta)\tau\big)^2}=4q(\tau),
\end{equation*}
and we arrive at
\begin{equation*}
i\partial_\tau Z+\Delta Z+|Z|^{p-1}Z+i\varepsilon\Lambda Z=\e^{-ib|y|^2}\Big[i\partial_\tau\widetilde Z+\Delta\widetilde Z+|\widetilde Z|^{p-1}\widetilde Z+q(\tau)|y|^2\widetilde Z\Big].
\end{equation*}
It remains to express the phase in the original variables. Since $\varepsilon\lambda^2=\lambda'/\lambda=-\frac{\beta}{1+t}$, we have $-b|\lambda x|^2=\frac{\beta|x|^2}{4(1+t)}$, which gives \eqref{eq:Ztilde}. The identity for the source term follows from the two displayed relations, and the conservation of the $L^2$ norm follows from the change of variables $y=\lambda x$.
\end{proof}

\begin{remark}\rm \label{rem:reduction}
(i) For $p<p_{\rm s}(a)$, that is, $\beta<\frac12$, the reduced equation \eqref{eq:repulsive} is defined for all $\tau\ge0$, the potential decays like $\tau^{-2}$, and
\begin{equation}\label{eq:intq}
\int_0^\infty q(\tau)\dd\tau=\frac{\beta(1-\beta)}{4(1-2\beta)}<\infty .
\end{equation}
Since $\tau(t)\sim\frac{t^{1-2\beta}}{1-2\beta}$, the reduced time $\tau$ runs much more slowly than $t$, and a profile which is stationary in the variables $(\tau,y)$ has width $(1+t)^\beta=o(t)$ in the original variables.

(ii) At the endpoint $p=p_{\rm s}(a)$, we have $\beta=\frac12$, $\tau=\log(1+t)$ and $q\equiv\frac1{16}$. The reduced equation \eqref{eq:repulsive} is then the \emph{autonomous} focusing NLS equation with the repulsive harmonic potential $\frac1{16}|y|^2$, studied by Carles \cite{Carles2003}. This reflects the scale invariance of Lemma \ref{lem:scaling}.

(iii) For $p>p_{\rm s}(a)$, we have $\tau_\infty<\infty$ and, unless $\beta=1$, $|q(\tau)|$ blows up like $(\tau_\infty-\tau)^{-2}$ (the potential becomes attractive when $\beta>1$). The reduced equation then lives on a finite interval, which is consistent with the fact that all solutions scatter (Theorem \ref{thm:sub}(ii)).

(iv) The pseudo-conformal (lens) transform used in \cite{BamriTayachi2023,Carles2002,Tao2009}, namely $v(t,x)=(1+t)^{-N/2}\e^{i|x|^2/(4(1+t))}W\big(\frac t{1+t},\frac x{1+t}\big)$, leads instead to $i\partial_sW+\Delta W+(1-s)^{-(2-\sigma-\alpha)}|W|^{p-1}W=0$ on $s\in[0,1)$. The coupling is singular at $s=1$ in the whole range $p<1+\frac4{N+2a}$. The reduction \eqref{eq:Ztilde} corresponds to the dispersion rate $(1+t)^\beta$ instead of $1+t$, and it is adapted to the solitary waves discussed below.
\end{remark}

\subsection{Rigorous consequences}
The first result shows that scattering forces a solution to spread at a scale comparable to $t$. It holds for every damping and every $p$.

\begin{proposition}\label{prop:concentration}
Let $v\in C([0,\infty);\Lt)$ and assume that there exists $v_+\in\Lt$ such that $\|\e^{-it\Delta}v(t)-v_+\|_{L^2}\to0$. Then, for every function $\rho:(0,\infty)\to(0,\infty)$ with $\rho(t)/t\to0$,
\begin{equation}\label{eq:noconc}
\lim_{t\to\infty}\ \sup_{x_0\in\R^N}\|v(t)\|_{L^2(B(x_0,\rho(t)))}=0 .
\end{equation}
In particular, if $u$ is a solution of \eqref{general} such that $\e^{A(t)}\|u(t)\|_{L^2(B(x(t),\rho(t)))}\ge c>0$ along a sequence of times $t\to\infty$, for some $x(t)\in\R^N$ and some $\rho(t)=o(t)$, then $u$ does not scatter in $L^2$. Moreover, if $1<p<1+\frac4N$ and the solution $u$ of \eqref{general} with $u_0\neq0$ scatters in $L^2$, then $E(u(t))>0$ for all sufficiently large $t$.
\end{proposition}

\begin{proof}
Let $Z(t)=\Mo(t)\Do(t)\hat v_+$. As in Section \ref{S6}, we have $\|v(t)-Z(t)\|_{L^2}\le\|\e^{-it\Delta}v(t)-v_+\|_{L^2}+\|(\Mo(t)-1)v_+\|_{L^2}\to0$. Moreover, the change of variables $x=2t\xi$ gives
\begin{equation*}
\|Z(t)\|_{L^2(B(x_0,\rho))}=\|\hat v_+\|_{L^2(B(\frac{x_0}{2t},\frac{\rho}{2t}))}\le\sup_{\xi_0\in\R^N}\|\hat v_+\|_{L^2(B(\xi_0,\frac{\rho}{2t}))}.
\end{equation*}
For every $g\in\Lt$, the quantity $\sup_{\xi_0}\|g\|_{L^2(B(\xi_0,r))}$ tends to zero as $r\to0$. Indeed, given $\eta>0$, we choose $g_\eta\in C_c(\R^N)$ with $\|g-g_\eta\|_{L^2}\le\eta$, and we get $\|g\|_{L^2(B(\xi_0,r))}\le\eta+\|g_\eta\|_{L^\infty}|B(0,r)|^{1/2}$. Applying this with $g=\hat v_+$ and $r=\rho(t)/(2t)\to0$, we obtain \eqref{eq:noconc}. The second assertion follows from \eqref{eq:noconc} applied to $v=\e^{A}u$.

For the last assertion, let $u$ scatter in $L^2$ with $u_0\neq0$, and let $v_+$ be the scattering state of $v=\e^{A}u$. By Proposition \ref{prop:apriori-sub}, $V:=\sup_t\|v(t)\|_{H^1}<\infty$. Hence $\e^{-it\Delta}v(t)$ is bounded in $\Ha$ and converges to $v_+$ in $\Lt$. We infer that $v_+\in\Ha$ and that $\nabla\big(\e^{-it\Delta}v(t)\big)\rightharpoonup\nabla v_+$ weakly in $\Lt$. In particular, $\liminf_{t\to\infty}\|\nabla v(t)\|_{L^2}\ge\|\nabla v_+\|_{L^2}>0$, since $v_+\neq0$, and $\|\nabla v_+\|_{L^2}\le V$. Furthermore, \eqref{Nirenberg} gives
\begin{equation*}
\|v(t)-\e^{it\Delta}v_+\|_{L^{p+1}}^{p+1}\le K\|v(t)-\e^{it\Delta}v_+\|_{L^2}^{p+1-\sigma}(2V)^{\sigma}\to0,
\end{equation*}
and $\|\e^{it\Delta}v_+\|_{L^{p+1}}\to0$ by \eqref{dispersive}, the density of $\cS(\R^N)$ in $\Ha$ and the embedding $\Ha\hookrightarrow L^{p+1}$. Hence $\|v(t)\|_{L^{p+1}}\to0$. Since $h\le1$, we conclude that $\liminf_{t\to\infty}\mathcal H(t)\ge\|\nabla v_+\|_{L^2}^2>0$, and $E(u(t))=\e^{-2A(t)}\mathcal H(t)$ is positive for large $t$.
\end{proof}

The second result is a pseudo-conformal law for \eqref{V}. It shows that, in a range which contains the whole gap, every solution with data in $\Sigma$ spreads at least at the rate $(1+t)^\beta$.

\begin{proposition}\label{prop:pseudoconf}
Let $\theta=1$, $a>0$, $(p-1)(N+2a)\le4$, and let $\beta$ be given by \eqref{eq:alphabeta}. Let $u_0\in\Sigma$, and let $v=\e^{A}u$, where $u$ is the global solution of \eqref{Main-New}. Set $J(t)=x+2it\nabla$ and
\begin{equation*}
\mathcal G(t)=\|J(t)v(t)\|_{L^2}^2-\frac{8t^2h(t)}{p+1}\|v(t)\|_{L^{p+1}}^{p+1}.
\end{equation*}
Then $\mathcal G$ is nonincreasing on $[0,\infty)$, and $\mathcal G(0)=\|xu_0\|_{L^2}^2$. Consequently,
\begin{equation}\label{eq:Jbound}
\|J(t)v(t)\|_{L^2}\le\max\Big\{\sqrt2\,\|xu_0\|_{L^2},\ (2B)^{\frac1{2-\sigma}}(1+t)^{1-\beta}\Big\},\qquad B=\frac{2^{3-\sigma}K}{p+1}m^{p+1-\sigma},
\end{equation}
and, for every $r\in[2,2^*]$ with $r<\infty$ if $N=2$, there exists $C$ depending on $N,p,a,r,m$ and $\|xu_0\|_{L^2}$ such that, with $\delta=N\big(\frac12-\frac1r\big)$,
\begin{equation}\label{eq:pcdecay}
\|v(t)\|_{L^r}\le C(1+t)^{-\beta\delta},\qquad \sup_{x_0\in\R^N}\|v(t)\|_{L^2(B(x_0,\rho))}\le C\Big(\frac{\rho}{(1+t)^{\beta}}\Big)^{\delta},\qquad t\ge1,\ \rho>0 .
\end{equation}
\end{proposition}

\begin{proof}
Note first that $(p-1)(N+2a)\le4$ and $a>0$ imply $p<1+\frac4N$, so that $\sigma<2$ and $u$ is global by Theorem \ref{thm:sub}(i). Note also that, by \eqref{eq:equivalences}, the assumption $(p-1)(N+2a)\le4$ means $\beta\le1$. For $t>0$, set $w(t)=\e^{-i|x|^2/(4t)}v(t)$. Then $|w|=|v|$ and $J(t)v=2it\,\e^{i|x|^2/(4t)}\nabla w$, so that $\|J(t)v\|_{L^2}=2t\|\nabla w\|_{L^2}$. Using $\nabla\e^{i|x|^2/(4t)}=\frac{ix}{2t}\e^{i|x|^2/(4t)}$, the computation of Step 2 in the proof of Proposition \ref{prop:reduction}, with $\varepsilon=0$ and $b=-\frac1{4t}$, shows that $w$ solves
\begin{equation*}
i\partial_tw+\Delta w+h|w|^{p-1}w+\frac it\Lambda w=0 .
\end{equation*}
(Indeed, with $\varepsilon=0$, the first-order terms combine into $-4ib\Lambda w=\frac it\Lambda w$, while the coefficient of $|x|^2w$ is $\partial_tb-4b^2=\frac1{4t^2}-\frac1{4t^2}=0$.) Let $\mathcal H_w=\|\nabla w\|_{L^2}^2-\frac{2h}{p+1}\|w\|_{L^{p+1}}^{p+1}$, so that $\mathcal G=4t^2\mathcal H_w$. Since $\Lambda$ is the generator of $L^2$-dilations, we have ${\rm Re}(-\Delta w,\Lambda w)=\|\nabla w\|_{L^2}^2$ and ${\rm Re}(|w|^{p-1}w,\Lambda w)=\frac{\sigma}{p+1}\|w\|_{L^{p+1}}^{p+1}$. Arguing as in the proof of Proposition \ref{prop:identities}, we obtain
\begin{equation*}
\frac{\dd}{\dd t}\mathcal H_w=-\frac2t\Big(\|\nabla w\|_{L^2}^2-\frac{\sigma h}{p+1}\|w\|_{L^{p+1}}^{p+1}\Big)-\frac{2h'}{p+1}\|w\|_{L^{p+1}}^{p+1}.
\end{equation*}
Therefore,
\begin{equation*}
\mathcal G'(t)=8t\mathcal H_w+4t^2\mathcal H_w'=\frac{8t\,h(t)\|v(t)\|_{L^{p+1}}^{p+1}}{p+1}\Big(\sigma-2+\frac{(p-1)a\,t}{1+t}\Big),
\end{equation*}
where we used \eqref{h-prime}. Since $\frac{t}{1+t}<1$, the bracket is negative as soon as $\sigma-2+\alpha\le0$, that is, $(p-1)(N+2a)\le4$. The computation is justified for $\Sigma$ solutions by the regularisation argument of \cite[Section 7.2]{Cazenave2003}, using the propagation of the $\Sigma$-regularity for \eqref{V}, which holds as in \cite[Section 6.5]{Cazenave2003} because $h$ is bounded. As $t\downarrow0$, we have $J(t)v(t)\to xu_0$ in $\Lt$ and $t^2h(t)\|v(t)\|_{L^{p+1}}^{p+1}\to0$, so that $\mathcal G(0)=\|xu_0\|_{L^2}^2$.

Let $Y(t)=\|J(t)v(t)\|_{L^2}^2$. By \eqref{Nirenberg} applied to $w$, we have $\|v\|_{L^{p+1}}^{p+1}\le Km^{p+1-\sigma}(2t)^{-\sigma}Y^{\sigma/2}$. Since $\mathcal G(t)\le\|xu_0\|_{L^2}^2$, we deduce that $Y\le\|xu_0\|_{L^2}^2+Bt^{2-\sigma}h(t)Y^{\sigma/2}$. As in the proof of Proposition \ref{prop:apriori-sub}, applied at each fixed $t>0$, this gives $Y\le\max\{2\|xu_0\|_{L^2}^2,(2Bt^{2-\sigma}h(t))^{2/(2-\sigma)}\}$. Since $(t^{2-\sigma}h)^{2/(2-\sigma)}=t^2(1+t)^{-2\beta}\le(1+t)^{2-2\beta}$, we obtain \eqref{eq:Jbound}. Finally, the Gagliardo--Nirenberg inequality $\|f\|_{L^r}\le C\|f\|_{L^2}^{1-\delta}\|\nabla f\|_{L^2}^{\delta}$, applied to $w$, gives $\|v(t)\|_{L^r}\le Cm^{1-\delta}(2t)^{-\delta}\|J(t)v(t)\|_{L^2}^{\delta}$. For $t\ge1$, we have $(2t)^{-\delta}\le(1+t)^{-\delta}$, and \eqref{eq:Jbound} together with $\beta\le1$ yields the first estimate in \eqref{eq:pcdecay}. The second one follows from H\"older's inequality, $\|v\|_{L^2(B(x_0,\rho))}\le|B(0,\rho)|^{\frac12-\frac1r}\|v\|_{L^r}=|B(0,1)|^{\frac\delta N}\rho^{\delta}\|v\|_{L^r}$.
\end{proof}

Propositions \ref{prop:concentration} and \ref{prop:pseudoconf} delimit the possible obstructions to scattering of concentration type in the gap. Scattering is incompatible with concentration at scales $o(t)$, while no solution with data in $\Sigma$ can concentrate at scales $o(t^\beta)$. Consequently, for data in $\Sigma$, an obstruction of this type can only come from a solution keeping a fixed amount of mass at scales between $t^\beta$ and $o(t)$. We stress that Proposition \ref{prop:concentration} only gives a sufficient condition for non-scattering; we do not claim that every non-scattering solution concentrates. The solitary waves considered below concentrate exactly at the smallest admissible scale $t^\beta$.

\subsection{Adiabatic solitary waves}
Let $R$ be the unique positive radial solution of $-\Delta R+R=R^p$ \cite{Kwong1989}, and set $R_\omega(y)=\omega^{\frac1{p-1}}R(\sqrt\omega\,y)$ for $\omega>0$. Then $-\Delta R_\omega+\omega R_\omega=R_\omega^p$ and $\|R_\omega\|_{L^2}^2=\omega^{\frac{2-\sigma}{p-1}}\|R\|_{L^2}^2$, so that every mass is attained. The function $\e^{i\omega\tau}R_\omega$ solves \eqref{eq:repulsive} without the potential term. Going back to the original variables through \eqref{eq:Ztilde}, we define
\begin{equation}\label{eq:vapp}
v_{\rm app}(t,x)=\e^{i\omega\tau(t)}\,(1+t)^{-\frac{N\beta}2}\exp\Big(\frac{i\beta|x|^2}{4(1+t)}\Big)R_\omega\Big(\frac{x}{(1+t)^{\beta}}\Big),\qquad u_{\rm app}=\e^{-A}v_{\rm app}=(1+t)^{-a}v_{\rm app}.
\end{equation}

\begin{proposition}\label{prop:approx}
Let $\theta=1$, $a>0$, $1<p<1+\frac4N$ and $\omega>0$. Then $\|v_{\rm app}(t)\|_{L^2}=\|R_\omega\|_{L^2}$, and $v_{\rm app}$ solves
\begin{equation*}
i\partial_tv_{\rm app}+\Delta v_{\rm app}+h(t)|v_{\rm app}|^{p-1}v_{\rm app}=\mathfrak r(t),\qquad \|\mathfrak r(t)\|_{L^2}=\frac{q(\tau(t))}{(1+t)^{2\beta}}\,\big\||y|^2R_\omega\big\|_{L^2}.
\end{equation*}
Moreover, if $p\neq1+\frac4{N+2a}$, then $\int_0^\infty\|\mathfrak r(t)\|_{L^2}\dd t<\infty$ if and only if $p<p_{\rm s}(a)$. In that case, for every $T\ge0$,
\begin{equation}\label{eq:errortail}
\int_T^\infty\|\mathfrak r(t)\|_{L^2}\dd t=\frac{\beta(1-\beta)\,\||y|^2R_\omega\|_{L^2}}{4(1-2\beta)}\,(1+T)^{-(1-2\beta)} .
\end{equation}
Finally, if $p<1+\frac4{N+2a}$, there is no $v_+\in\Lt$ such that $\|\e^{-it\Delta}v_{\rm app}(t)-v_+\|_{L^2}\to0$.
\end{proposition}

\begin{proof}
The function $\widetilde Z=\e^{i\omega\tau}R_\omega$ solves \eqref{eq:repulsive} with the source term $F=q(\tau)|y|^2\widetilde Z$. By Proposition \ref{prop:reduction}, $v_{\rm app}$ solves \eqref{V} with the source $\mathfrak r(t,x)=\lambda^{\frac N2+2}\exp\big(\frac{i\beta|x|^2}{4(1+t)}\big)F(\tau,\lambda x)$, and $\|\mathfrak r(t)\|_{L^2}=\lambda^2\|F(\tau)\|_{L^2}=\lambda(t)^2q(\tau(t))\||y|^2R_\omega\|_{L^2}$. Since $\dd\tau=\lambda^2\dd t$, we get
\begin{equation*}
\int_T^\infty\|\mathfrak r(t)\|_{L^2}\dd t=\||y|^2R_\omega\|_{L^2}\int_{\tau(T)}^{\tau_\infty}q(\tau)\dd\tau .
\end{equation*}
If $\beta<\frac12$, then $\tau_\infty=\infty$ and $\int_{\tau(T)}^\infty q\dd\tau=\frac{\beta(1-\beta)}{4(1-2\beta)}\big(1+(1-2\beta)\tau(T)\big)^{-1}$, which gives \eqref{eq:errortail}. If $\beta=\frac12$, then $q\equiv\frac1{16}$ and $\tau_\infty=\infty$, so that the integral diverges. If $\beta>\frac12$, then $q(\tau)=\frac{\beta(1-\beta)}{4(2\beta-1)^2(\tau_\infty-\tau)^2}$ is not integrable near $\tau_\infty$, because $\beta(1-\beta)\neq0$ for $\beta\neq1$. By \eqref{eq:equivalences}, the value $\beta=1$ corresponds to $p=1+\frac4{N+2a}$, which is excluded, and the integral is finite exactly when $p<p_{\rm s}(a)$. Finally, for every $\rho>0$, the mass of $v_{\rm app}(t)$ in the ball $B(0,\rho(1+t)^\beta)$ equals $\|R_\omega\|_{L^2(B(0,\rho))}$, which is positive and independent of $t$. If $p<1+\frac4{N+2a}$, then $\beta<1$ by \eqref{eq:equivalences}, so that $\rho(1+t)^\beta=o(t)$, and Proposition \ref{prop:concentration} shows that $v_{\rm app}$ is not asymptotically free.
\end{proof}

\begin{remark}\rm 
The exponent $p=1+\frac4{N+2a}$ excluded in Proposition \ref{prop:approx} corresponds to $\beta=1$. There, $q\equiv0$, the function $v_{\rm app}$ is an exact solution, and \eqref{eq:Ztilde} is the classical pseudo-conformal transform, since $\lambda(t)=(1+t)^{-1}$ and $\tau(t)=\frac t{1+t}$. In that case, $v_{\rm app}$ is the pseudo-conformal image of a solitary wave, and it spreads at the rate $t$, as expected.
\end{remark}

Proposition \ref{prop:approx} is not an existence result. An $L^1_tL^2_x$-small error does not by itself produce an exact solution nearby, because $v_{\rm app}$ has infinite Strichartz norms on $[0,\infty)$, so that the perturbative argument of Section \ref{S5} does not apply. Indeed, $\|v_{\rm app}(t)\|_{L^r}=(1+t)^{-\beta\delta}\|R_\omega\|_{L^r}$ with $\delta=N(\frac12-\frac1r)=\frac2q$, and $t\mapsto(1+t)^{-2\beta/q}$ is not in $L^q(0,\infty)$ when $\beta\le\frac12$. The relevant stability theory is orbital rather than perturbative. In the variables of Proposition \ref{prop:reduction}, the question becomes the following: does the orbital stability of the ground state $\e^{i\omega\tau}R_\omega$ \cite{CazenaveLions1982,Weinstein1985,Weinstein1986} survive the slowly decaying repulsive potential $q(\tau)|y|^2$? The dynamics of solitary waves in slowly varying external potentials is well understood in related settings \cite{FGJS2004,HolmerZworski2008}, and the integrability \eqref{eq:intq} makes the perturbation summable in time. This leads to the following conjecture.

\begin{conjecture}[Adiabatic solitary waves]\label{conj:adiabatic}
Let $\theta=1$, $a>0$, $N\ge1$ and $1<p<p_{\rm s}(a)$. For every $\omega>0$, there exists a solution $u$ of \eqref{Main-New}, with radial data if $N\ge2$, such that, with $v=\e^{A}u$,
\begin{equation*}
\lim_{t\to\infty}\ \inf_{\gamma\in\R}\big\|v(t)-\e^{i\gamma}v_{\rm app}(t)\big\|_{L^2}=0 .
\end{equation*}
In particular, $u$ does not scatter in $L^2$. Consequently, for $p_{\rm ns}(a)<p<p_{\rm s}(a)$, scattering and non-scattering solutions coexist, and $p_{\rm s}(a)$ is the sharp threshold for the scattering of all solutions.
\end{conjecture}

The non-scattering part of the conclusion follows from the first part by Proposition \ref{prop:approx} and Proposition \ref{prop:concentration}. The conjecture is consistent with all the results of this paper. Indeed, the solitary waves of Conjecture \ref{conj:adiabatic} have width $(1+t)^\beta$, which is the smallest scale allowed by Proposition \ref{prop:pseudoconf}. Moreover, a direct computation gives
\begin{equation*}
\begin{aligned}
(1+t)^{2\beta}\Big(\|\nabla v_{\rm app}\|_{L^2}^2-\frac{2h}{p+1}\|v_{\rm app}\|_{L^{p+1}}^{p+1}\Big)&=E(R_\omega)+\frac{\varepsilon^2}4\|yR_\omega\|_{L^2}^2\\
&\xrightarrow[t\to\infty]{}E(R_\omega)=-\frac{2-\sigma}{p+1}\|R_\omega\|_{L^{p+1}}^{p+1}<0,
\end{aligned}
\end{equation*}
where we used the Pohozaev identity $\|\nabla R_\omega\|_{L^2}^2=\frac{\sigma}{p+1}\|R_\omega\|_{L^{p+1}}^{p+1}$, and where $\varepsilon=-\beta(1+t)^{2\beta-1}\to0$ is the function of Proposition \ref{prop:reduction}. The conjectured solutions should therefore have negative energy for large times, which is compatible with the last assertion of Proposition \ref{prop:concentration}. The mechanism is the classical adiabatic one (see \cite{KivsharMalomed1989} for the formal adiabatic theory of perturbed solitons). The internal frequency of the solitary wave is $\omega\lambda(t)^2=\omega(1+t)^{-2\beta}$, while the coupling varies on the time scale $1+t$. The number of internal oscillations, $\omega\tau(t)$, tends to infinity if and only if $\beta\le\frac12$, and the accumulated error is finite if and only if $\beta<\frac12$.

\begin{figure}[ht]
\centering
\begin{tikzpicture}
\begin{axis}[width=11.5cm,height=6.2cm,xmin=0,xmax=400,ymin=-0.02,ymax=1.05,
 xlabel={reduced time $\tau$},ylabel={$\|\widetilde Z(\tau)\|_{L^2(|y|<15)}^2/\|\widetilde Z(0)\|_{L^2}^2$},
 legend style={at={(0.98,0.55)},anchor=east,font=\footnotesize,fill=white},
 tick label style={font=\footnotesize},label style={font=\small}]
\addplot[thick,green!50!black] coordinates {(0,1.0000) (10,0.9996) (20,0.9987) (30,0.9983) (40,0.9981) (50,0.9981) (60,0.9980) (80,0.9979) (100,0.9979) (150,0.9979) (200,0.9978) (250,0.9978) (300,0.9978) (350,0.9978) (400,0.9978)};
\addlegendentry{$a=0.1$, $\beta=0.2$ (inside the gap)}
\addplot[thick,dashed,orange!85!black] coordinates {(0,1.0000) (2,1.0000) (4,0.9999) (6,0.9927) (8,0.9752) (10,0.9586) (12,0.9482) (14,0.9434) (16,0.9420) (18,0.9418) (20,0.9417) (30,0.9416) (40,0.9416) (100,0.9416) (200,0.9416) (300,0.9416) (400,0.9416)};
\addlegendentry{$a=0.225$, $\beta=0.45$ (inside the gap)}
\addplot[thick,dotted,red!75!black] coordinates {(0,1.0000) (2,1.0000) (4,0.9991) (6,0.9684) (8,0.8750) (10,0.6887) (12,0.4003) (14,0.1691) (16,0.0637) (18,0.0235) (20,0.0087) (22,0.0032) (24,0.0012) (26,0.0004) (30,0.0001) (40,0.0000) (100,0.0000) (200,0.0000) (400,0.0000)};
\addlegendentry{$a=0.25$, $\beta=\frac12$ (endpoint $p=p_{\rm s}(a)$)}
\end{axis}
\end{tikzpicture}
\caption{Numerical solution of the reduced equation \eqref{eq:repulsive} for $N=1$, $p=3$ and $\widetilde Z(0)=R_1=\sqrt2\,{\rm sech}$, for three values of $a$. Plotted is the fraction of the mass contained in $\{|y|<15\}$. For $N=1$ and $p=3$, one has $\beta=2a$, $p_{\rm ns}(a)<3$ for all $a>0$, and $3<p_{\rm s}(a)$ if and only if $a<\frac14$. Inside the gap, the solitary wave remains trapped, after shedding a small amount of radiation. At the endpoint, it is destroyed. The reduced time $\tau=400$ corresponds to $t\approx9.3\cdot10^3$ for $\beta=0.2$, to $t\approx1.3\cdot10^{16}$ for $\beta=0.45$, and to $t=\e^{400}-1$ for $\beta=\frac12$.}
\label{fig:numerics}
\end{figure}
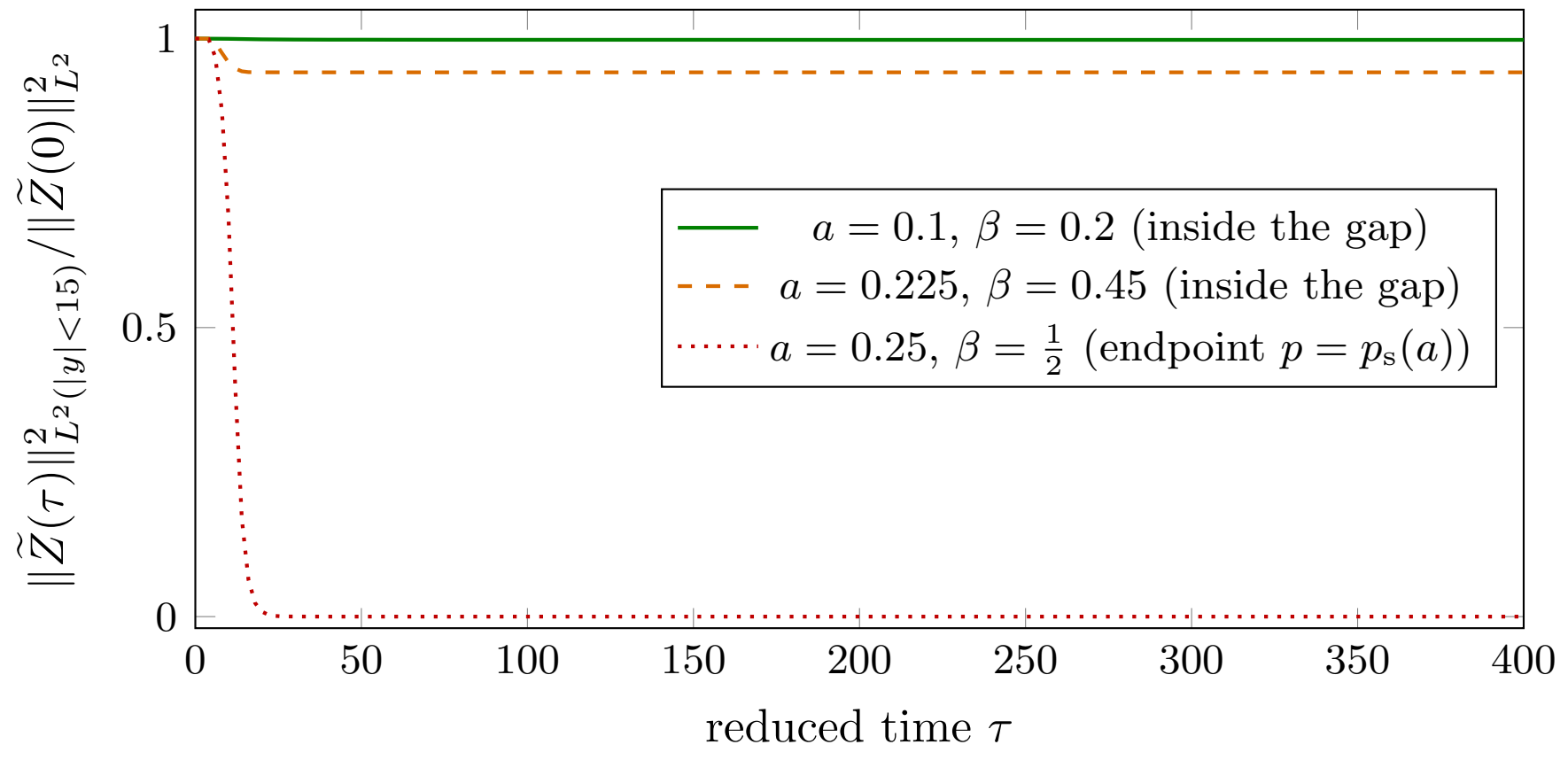

\begin{remark}\rm [Numerical evidence]\label{rem:numerics}
We solved \eqref{eq:repulsive} for $N=1$ and $p=3$ by a second-order split-step Fourier method on $[-300,300]$, with $2^{14}$ grid points, time step $\Delta\tau=5\cdot10^{-3}$ and an absorbing layer on the outer $15\%$ of the domain. For $N=1$ and $p=3$, every $a\in(0,\frac14)$ lies in the open gap, and Theorem \ref{thm:sub}(ii) applies exactly when $a>\frac14$. Figure \ref{fig:numerics} shows that, for $a=0.1$ and $a=0.225$, the ground state $R_1$ remains trapped, with a local mass fraction of $0.998$ and $0.942$ at $\tau=400$, respectively, and that this persists up to $\tau=1.2\cdot10^3$ (equivalently $t\approx5.8\cdot10^4$ for $a=0.1$). Halving the time step changes these values by less than $10^{-5}$. In longer runs, a small symmetry-breaking perturbation produced by rounding errors makes the solitary wave drift at the rate $\tau^{(1-\beta)/(1-2\beta)}$ in the variable $y$, without any loss of amplitude. This rate corresponds exactly to a constant velocity in the original variable $x$, that is, to a harmless Galilean motion. At the endpoint $a=\frac14$, the solitary wave of mass $\|R_1\|_{L^2}^2=4$ is destroyed by $\tau\approx20$. Solitary waves with larger mass ($\omega=4$ or $9$) survive much longer, which suggests a slow, tunnelling-type decay. For this reason, we leave the endpoint open. Finally, the profile $\frac12R_1$, whose mass is $1$, disperses for $a=0.1$, in agreement with small data scattering.
\end{remark}

\subsection{Comparison with small data scattering}
For comparison, consider the small data theory of \cite[Corollary 1.6]{BamriTayachi2023} ($N\ge2$) with $\gamma=2a$. There, scattering holds for $p-1>\alpha_1(2a)$ (within the ranges of $p$ specified in \cite[Corollary 1.6]{BamriTayachi2023}, all of which satisfy $p-1<\frac4{N-2}$), for data whose free evolution is small in a Lorentz--Strichartz norm, where $\alpha_1(\gamma)$ is the positive root of
\begin{equation}\label{eq:quadratic}
(N+\gamma)\alpha^2+(N-2+2\gamma)\alpha-4=0 .
\end{equation}
For $\gamma=0$, $1+\alpha_1(0)$ is the classical Strauss exponent. A direct computation shows that $\frac2{N+2a}<\alpha_1(2a)<\frac4{N+4a}$, so that the small data threshold lies strictly inside the gap. Indeed, for $N\ge2$, the left-hand side of \eqref{eq:quadratic} with $\gamma=2a$ is increasing in $\alpha>0$, since both coefficients $N+2a$ and $N-2+4a$ are nonnegative. Evaluated at $\alpha=\frac4{N+4a}$, it equals $\frac{8N}{(N+4a)^2}>0$, and evaluated at $\alpha=\frac2{N+2a}$, it equals $-\frac{2N}{N+2a}<0$. Thus, for small data in the Lorentz--Strichartz sense, the range of scattering is larger than the one of Theorem \ref{thm:sub}(ii), whereas Theorem \ref{thm:sub}(ii) applies to data of arbitrary size. This smallness condition is satisfied by sufficiently small data in $\Sigma$: this is shown in \cite[Remark 1.3(4)]{BamriTayachi2023} for the setting of \cite[Theorem 1.2]{BamriTayachi2023}, with $q\ge d$ in their notation and under the condition $p-1>\alpha_0(\delta)$, where $\alpha_0(\delta)$ is the positive root of $N\alpha^2+(N-2+2\delta)\alpha-4(1-\delta)=0$. A direct computation shows that the lower bound imposed on $\delta$ in \cite[Corollary 1.6]{BamriTayachi2023}, namely $\delta>\frac{4-(N-2)\alpha-N\alpha^2}{2(\alpha+2)}$ with $\alpha=p-1$, is equivalent to $p-1>\alpha_0(\delta)$, and that corollary allows any $q\in[1,\infty)$. Hence, for $N\ge2$, small data in $\Sigma$ already scatter for $p-1>\alpha_1(2a)$ within these ranges. Whether the same holds for small data in $\Ha$ without weight, or for $N=1$, remains open.

\subsection{Proof of Theorem \ref{thm:gap}}
Part (i) is Lemma \ref{lem:scaling}, together with \eqref{eq:equivalences}. Part (ii) is Proposition \ref{prop:reduction}, together with Remark \ref{rem:reduction}(i)--(ii). The first assertion of part (iii) is the second assertion of Proposition \ref{prop:concentration}. Its second assertion is the second estimate in \eqref{eq:pcdecay}, since $v=\e^{A}u$, with $\delta=N(\frac12-\frac1r)$ for any fixed $r\in(2,2^*]$, $r<\infty$ if $N=2$. Finally, part (iv) is Proposition \ref{prop:approx}, the function $v_{\rm app}$ being given by \eqref{eq:vapp}.

\end{document}